%% file: main.tex
\documentclass[11pt, a4paper, oneside,reqno]{amsart}

\input{AMIN_style}

\graphicspath{{Fig/}}

\usepackage[square,numbers]{natbib}
\usepackage{bm}

\newcommand{\bo}{T}

\DeclareMathOperator{\id}{Id}
\newcommand{\sbo}{\widehat{T}}

\newcommand{\hf}{\widehat{f}}
\newcommand{\gr}{\nabla}

\newcommand{\e}{\boldsymbol{1}}

\title[]{Control and Reinforcement Learning through the Lens of Optimization: An Algorithmic Perspective}

\author[]{Tolga~Ok$^{1}$, Arman~Sharifi~Kolarijani$^1$, Mohamad Amin Sharifi Kolarijani$^1$, and Peyman~Mohajerin~Esfahani$^{1,2}$ \\
\\
$^1$Delft Center for Systems and Control, Delft University of Technology, The Netherlands\\
$^2$Department of Mechanical and Industrial Engineering, University of Toronto, Canada
}
\thanks{Correspondence to: Tolga Ok $<$\texttt{t.ok@tudelft.nl}$>$.}
\thanks{This work was partially supported by the European Research Council (ERC) project TRUST-949796, the Horizon Europe Pathfinder Open project RELIEVE-101099481, and the NSERC Discovery grant RGPIN-2025-06544.} 

\begin{document} 

\begin{abstract}
The connection between control algorithms for Markov decision processes and optimization algorithms has been implicitly and explicitly exploited since the introduction of dynamic programming algorithm by Bellman in the 1950s. 
Recently, this connection has attracted a lot of attention for developing new control algorithms inspired by well-established optimization algorithms. 
In this paper, we make this analogy explicit across four problem classes with a unified solution characterization. 
This novel framework, in turn, allows for a systematic transformation of algorithms from one domain to the other. 
In particular, we identify equivalent optimization and control algorithms that have already been pointed out in the existing literature, but mostly in a scattered way. 
We also discuss the issues arising in providing theoretical convergence guarantees for these new control algorithms and provide simple yet effective techniques to solve them. 
The provided framework and techniques then lay out a concrete methodology for developing new convergent control algorithms. 

\smallskip
\smallskip
\noindent\textsc{Keywords:} Markov decision processes, stochastic optimal control, reinforcement learning, optimization, accelerated algorithms, quasi-Newton methods 

\end{abstract}

\maketitle

\section{Introduction}\label{sec:intro}

Markov decision processes (MDPs) have become a standard mathematical framework for the formulation of stochastic optimal control problems, that is, the problem of dynamic decision-making under uncertainty. 
This popularity can be to a large extent attributed to the Bellman principle of optimality and the dynamic programming algorithm~\cite{bellman1957markovian}. 
In particular, the fixed-point characterization of the optimal value function of an MDP has led to the development of a large class of iterative value-based algorithms, such as value iteration and policy iteration. 
This fixed-point characterization has also been used for identifying the fundamental connection between algorithms for optimal control of MDPs and those for optimization. 
A classic example is the policy iteration algorithm for MDPs which is an instance of the Newton method~\cite{puterman1979convergence, bertsekas2022lessons}. 
Recently, the connection between optimization and control problems has attracted a lot of attention for developing new control algorithms, with faster convergence and/or lower complexity, inspired by their counterparts for solving optimization problems. 
For instance, the modifications to the value iteration algorithm in~\cite{goyal2019first} and the Q-learning algorithm in~\cite{weng2020momentum} are inspired by momentum-based acceleration in optimization. 

The implicit connection between optimization algorithms and control algorithms for MDPs with a finite state-action space has also been studied more systematically. 
In~\cite{vieillard2019connections}, the authors look at the connection between \emph{constrained} convex optimization algorithms and control algorithms such as Frank-Wolfe algorithm~\cite{frank1956algorithm} and conservative policy iteration~\cite{kakade2002approximately}. A detailed comparison between \emph{deterministic} optimization algorithms and \emph{model-based}\footnote{In this paper, the terminologies of ``model-free'' and ``model-based'' indicate the {\em available information (oracle)}, i.e., whether we have access to the model or only the system trajectory (samples); see Section~\ref{sec:control problem} for more details. We note that this is different from the common terminologies in the RL literature where these terms refer to the {\em solution approach}, i.e., whether we identify the model along the way (model-based RL) or directly solve the Bellman equation to find the value function (model-free RL).} control algorithm is also provided in \cite{grand2021convex}, where the author looks at a wide range of optimization algorithms including gradient descent, accelerated gradient descent, Newton method, and quasi-Newton method and their counterparts for solving control problems. 

Motivated by these observations, this study provides an explicit framework unifying the tight link between optimization (stochastic and deterministic, respectively) and optimal control (model-based dynamic programming and model-free reinforcement learning, respectively). 
Specifically, this goal is achieved by exploiting the (expected) root-finding characterization of optimization problems and the (expected) fixed-point characterization of control problems. 
The framework yields an explicit transformation of deterministic (stochastic) convex optimization problems to model-based (model-free) control problems, and vice versa (Section~\ref{sec:problem}, Table~\ref{tab:equivalence}).
This explicit relationship, in turn, allows for a systematic transformation of algorithms from one domain to the other, which we will use to identify existing (and mostly known) equivalent algorithms for optimization and control (Section~\ref{sec:equiv_alg}, Table~\ref{tab:algorithms}). 
Unfortunately, the common formulation for control and optimization algorithms does not allow us to borrow existing theoretical results from one field to another.
Indeed, the similarity in the update rule does not carry over to the analysis, particularly when it comes to the convergence of these algorithms. 
Therefore, we also examine general tools for the convergence analysis of control algorithms, as well as some simple yet effective remedies for ensuring the convergence of novel control algorithms (Section~\ref{sec:convergence}, Theorems~\ref{thm:safeguardVI}, \ref{thm:backtrackVI}, \ref{thm:safeguardQL}). 
The provided framework and techniques then lay out a concrete methodology for developing new algorithms in one domain based on the existing algorithms in the other, along with theoretical convergence guarantees. 

In this paper, we restrict attention to control algorithms for \emph{finite} state–action MDPs.
For infinite (continuous) state–action spaces, aside from special cases like linear–quadratic regulators (LQR), practical computation typically requires finite-dimensional approximations.
One approach is to approximate at the model level by aggregating (discretizing) the state and action spaces, thereby reducing the problem to a finite MDP~\cite{ref:Bert_disc_75,ref:Powell_07}.
Another approach is to approximate the value function directly with a finite parameterization by minimizing (a proxy for) the residual of its Bellman fixed-point equation~\cite{ref:Bert_abstract, szepesvari2022algorithms}. Examples include linear parameterizations~\cite{ref:Bert_temporal, Tsitsiklis97}, nonlinear parameterizations using neural networks~\cite{ref:Bert_neuro-dynamic,SCHMIDHUBER2015,sutton2018reinforcement}, and max-plus methods~\cite{Bach20, GONCALVES2021109623, ref:FDP_TAC, ref:FDP_NeurIPS, McEn06, Liu2024fitted}.
There is also an alternative formulation in which the target function solves an infinite-dimensional linear program~\cite{ref:HL1}, enabling approximation via finite, tractable convex programs~\cite{ref:VanRoy_MP,ref:HL2,ref:Moh_SIOPT}.
Viewed this way, most of these approximation strategies can be framed as solving a finite-dimensional fixed-point problem or a convex optimization problem.

\noindent\textbf{Notations.} For a vector $v \in \R^n$, we use $v(i)$ and $[v](i)$ to denote its $i$-th element. 
Similarly, $M(i,j)$ and $[M](i,j)$ denote the element in row $i$ and column $j$ of the matrix $M \in \R^{m\times n}$. 
We use~$\cdot\tr$ to denote the transpose of a vector/matrix. 
$\norm{\cdot}_2$ and $\norm{\cdot}_{\infty}$ denote the 2-norm and $\infty$-norm of a vector, respectively. 
$\norm{\cdot}_2$ and $\norm{\cdot}_F$ denote the induced 2-norm and the Frobenius norm of a matrix, respectively. 
The identity operator is denoted by $\id$. 
We use $\e$ and $I$ to denote the all-ones vector and the identity matrix, respectively. 
We denote the $i$-th unit vector by $e_i$, that is, the vector with its $i$-th element equal to $1$ and all other elements equal to $0$. 

\section{Problems: Optimization vs. Control}
\label{sec:problem}
In this section, we provide the generic framework that connects optimization problems to control problems. 
In particular, we provide the explicit transformations between different variables and operations in these problems. 
Table~\ref{tab:equivalence} provides a condensed summary of this framework.

\input{Tab/equivalence_table_4}

\subsection{Optimization problem} 
We first look at the root-finding characterization of the solution to convex optimization problems. 
Consider the minimization problem 
\begin{equation}\label{eq:optimization_problem}
    \min_{x \in \R^{\ell}} \left\{ f(x) = \EE_{\xi}[\hf(x,\xi)] \right\},
\end{equation}
where $\hf : \R^{\ell} \times \Xi \to \R$ is a sample-wise function and $\xi$ is random vector with a fixed probability distribution $\PP$ over $\Xi$.
We assume that the resulting function $f : \R^{\ell} \to \R$ is twice continuously differentiable and strongly convex.
Much like the control problem, this problem can be considered in two settings:

\begin{itemize}
\item[(i)] 
    \textbf{Deterministic} optimization: Assume that $\PP$ in \eqref{eq:optimization_problem} is known and that the corresponding expectation can be computed. The unique minimizer~$x^\star$ is the root of the gradient operator, i.e.,  
    \begin{equation}\label{eq:fp_gradient}
        \gr f (x^\star) = 0.
    \end{equation}
    
\item[(ii)] 
    \textbf{Stochastic} optimization: Now assume that $\PP$ in \eqref{eq:optimization_problem} is unknown but can be sampled from. 
    In this case, the minimizer $x^\star$ satisfies the expected root-finding problem 
    \begin{equation}\label{eq:expected_fp_gradient}
       \EE_{\xi}[\gr \hf(x^\star,\xi)]  = 0, 
    \end{equation}
    where $\nabla$ now denotes the partial derivative w.r.t.~$x$. 
    We note that, above, there is an underlying assumption that the differentiation w.r.t.~$x$ and expectation w.r.t.~$\xi$ can be operated in any order.

\end{itemize}

\subsection{Control problem}\label{sec:control problem} 
The standard modeling framework for stochastic optimal control problems is the \emph{Markov decision process} (MDP). 
An MDP is characterized by the tuple $(\set{S},\set{A},\PP,c,\gamma)$, where
\begin{itemize}[label=$\bullet$, itemsep = 1mm, topsep = 0mm, leftmargin = 8mm]
    \item $\set{S}$ is the state space,
    \item $\set{A}$ is the action space,
    \item $\PP$ is the transition probability kernel such that $\PP(s^+|s,a)$ is the probability of the transition to state~$s^+$ given that the system is in state~$s$ and the chosen control is $a$ for each $(s,a,s^+) \in \set{S}\times\set{A}\times\set{S}$, 
    \item $c \in \R^{\set{S}\times\set{A}}$ is the stage cost function such that $c(s,a)$ is the cost of taking the control action~$a$ while the system is in state~$s$, and, 
    \item $\gamma \in (0,1)$ is the discount factor acting as a trade-off between short- and long-term costs.
\end{itemize}

Let us now fix a \emph{control policy}~$\pi:\set{S}\ra\set{A}$, i.e., a mapping from states to actions. 
The stage cost of the policy $\pi$ is denoted by $c^{\pi} \in \R^{\set{S}}$, where $c^{\pi}(s) \Let c\big(s,\pi(s)\big)$ for $s \in \set{S}$. 
The transition probability kernel of the resulting Markov chain under the policy $\pi$ is denoted by $\PP^{\pi}$, where $\PP^{\pi} (s^+|s) = \PP\big(s^+|s,\pi(s)\big)$ for $s,s^+\in \set{S}$. 
We also define $P^\pi \in [0,1]^{\set{S}\times\set{S}} $, with  $P^\pi(s,s^+) \Let \PP^{\pi}(s^+ | s)$ for each $(s,s^+) \in \set{S}\times\set{S}$, to be the corresponding transition probability \emph{matrix}. 
The value~$v^{\pi} \in \R^{\set{S}}$ of the policy~$\pi$ is then the expected, discounted, accumulative cost of following this policy over an infinite-horizon trajectory, that is,  
$$v^{\pi}(s) \Let \EE_{(s_{t})_{t=0}^\infty}\big[\ssum_{t=0}^\infty\gamma^t c(s_t,a_t)\mid s_0 = s,\ a_t = \pi(s_t)\big], \quad \forall s\in\set{S}.$$ 
Let us also define the \emph{action-value function} (a.k.a. the \emph{Q-function}) $q^{\pi} \in \R^{\set{S}\times\set{A}}$ for the policy~$\pi$ by 
$$q^{\pi}(s,a) \Let c(s,a) + \gamma\ \EE_{s^+} [v^{\pi}(s^+)\mid (s,a)],\quad  \forall(s,a)\in \set{S}\times\set{A},$$ 
so that $v^\pi(s) = q^\pi\big(s,\pi(s)\big)$ for each $s\in\set{S}$. Given a value function $v\in \R^{\set{S}}$, let us also define the \emph{greedy policy $\pi_v:\set{S}\ra\set{A}$ w.r.t.~$v$} denoted and given by 
$$\pi_v(s) \in \argmin_{a\in\set{A}}\left\{ c(s,a) + \gamma\ \EE_{s^+}\left[v(s^+)\mid(s,a)\right]\right\}, \quad \forall s\in\set{S}.$$ 
Similarly, for a Q-function~$q\in \R^{\set{S}\times\set{A}}$, we define the \emph{greedy policy $\pi_q:\set{S}\ra\set{A}$ w.r.t.~$q$} by 
$$\pi_q(s) \in \argmin_{a\in\set{A}} q(s,a),\quad  \forall(s,a)\in \set{S}\times\set{A}.$$ 
 
The problem of interest is to control the MDP optimally, that is, to find an optimal policy $\pi^*$ with the optimal (action-)value functions 
\begin{equation}\label{eq:optimal value func}
    v\opt = \min_{\pi} v^{\pi} \quad \text{and} \quad q\opt = \min_{\pi} q^{\pi},
\end{equation}
so that the expected, discounted, infinite-horizon cost is minimized. 
Let us also note that the optimal policy, that is a minimizer of the preceding optimization problems, is greedy w.r.t.~$v\opt$ and $q\opt$, i.e., $\pi\opt = \pi_{v\opt} = \pi_{q\opt}$. 

Interestingly, the optimal (action-)value functions introduced in~\eqref{eq:optimal value func} can be equivalently characterized as the fixed-point of two different operators each of which is useful depending on the available information (oracle): 

\begin{enumerate}[label = (\roman*), itemsep = 1mm, topsep = 0mm, leftmargin = 8mm]
\item 
\textbf{Model-based} control: When we have access to the transition kernel and the cost function, the problem is usually characterized by the fixed-point problem. 
    Indeed, defining the \emph{Bellman optimality operator}~$\bo:\R^{\set{S}}\ra\R^{\set{S}}$ by \begin{equation*} \label{eq:T}
        [\bo (v)](s) \Let \min_{a \in\set{A}} \left\{ c(s,a) + \gamma\ \EE_{s^+} \left[v(s^+)\mid(s,a)\right] \right\},\quad \forall s\in\set{S}, 
    \end{equation*}
    we have
    \begin{equation}\label{eq:fp_v}
        v\opt(s) = [\bo (v\opt)](s),\quad \forall s\in\set{S},
    \end{equation}
    or, equivalently, $v\opt = \bo(v\opt)$, that is, the optimal value function~$v\opt$ is the \emph{unique} fixed-point of~$\bo$. 
    The uniqueness follows from the fact that the operator $\bo$ is a $\gamma$-contraction in $\infty$-norm. 
    Observe that, in this case, $\bo$ can be exactly computed given the model of the underlying MDP. 
\item 
    \textbf{Model-free} control: In real applications, the model is often unknown, and instead, one can generate samples. 
    Examples of this are very large systems where identifying the model is prohibitively expensive but transitions between states can be observed and recorded, such as those in video games. 
    This problem has been studied extensively in the reinforcement learning community and is often formulated as the \emph{expected} fixed-point problem. 
    To provide this characterization, let us define the \emph{vector of sampled next states}~$\varsigma^+ \in \set{S}^{\set{S}\times \set{A}}$ where 
    \begin{equation}\label{eq:next_state_vector}
        \varsigma^+(s,a) \Let s^+ \sim \PP(\cdot|s,a),
    \end{equation}
    is a sample of the next state drawn from the distribution $\PP(\cdot|s,a)$ for each $(s,a)\in\set{S}\times\set{A}$. 
    Then, defining the \emph{sampled} Bellman optimality operator~$\sbo:\R^{\set{S}\times \set{A}}\times\set{S}^{\set{S}\times \set{A}}\ra\R^{\set{S}\times \set{A}}$ by\footnote{Strictly speaking, the provided sampled Bellman optimality operator is the empirical version of the Bellman optimality operator for the Q-function, given by $[\bar{\bo} (q)](s,a) \Let c(s,a) + \gamma \EE_{s^+} \left[\min_{a^+ \in\set{A}} q(s^+,a^+) \mid (s,a)\right]$ for each $(s,a)\in\set{S}\times\set{A}$.}
    \begin{equation*} \label{eq:rT}
        [\sbo (q,\varsigma^+)](s,a) \Let c(s,a) + \gamma \min_{a^+ \in \set{A}} q\big(\varsigma^+(s,a),a^+\big), \quad \forall (s,a)\in\set{S}\times\set{A}, 
    \end{equation*}
    we have
    \begin{equation}\label{eq:fp_q}
        q\opt(s,a) = \EE_{\varsigma^+(s,a)}\left[[\sbo (q\opt,\varsigma^+)](s,a) \right],\quad  \forall(s,a)\in \set{S}\times\set{A} ,
    \end{equation}
    or, equivalently, $q\opt = \EE_{\varsigma^+} [\sbo (q\opt,\varsigma^+)]$.

\end{enumerate}

\subsection{Transformation} 

In what follows, we focus on \emph{tabular MDPs with a finite state-action space}. We thus consider $\set{S} = \{1,2,\ldots,n \}$ and $\set{A} = \{1,2,\ldots,m \}$. 
Recall the Bellman optimality operator~$T : \R^{\set{S}}\ra\R^{\set{S}}$ with
\begin{align*}
    \bo (v) = \ssum_{s \in\set{S}} \ [\bo (v)](s) \cdot e_s,
\end{align*}
where $e_s \in \R^{\set{S}}$ is the unit vector for the state $s \in \set{S}$, 
and the sampled Bellman optimality operator~$\sbo : \R^{\set{S}\times\set{A}}\times\set{S}^{\set{S}\times\set{A}} \ra \R^{\set{S}\times\set{A}}$ with
\begin{align*}
    \sbo (q,\varsigma^+) = \ssum_{(s,a) \in\set{S}\times\set{A}} \ [\sbo (q,\varsigma^+)](s,a) \cdot e_{(s,a)},
\end{align*}
where $e_{(s,a)} \in \R^{\set{S}\times\set{A}}$ is the unit vector for the state-action pair~$(s,a) \in \set{S}\times\set{A}$, 
and $\varsigma^+ \in \set{S}^{\set{S}\times\set{A}}$ is the vector of sampled next states defined in~\eqref{eq:next_state_vector}. 

Before providing the equivalence relations between optimization and control problems, let us provide an important result for the Bellman optimality operator. 
\begin{Lem}[Jacobian of $\bo$]\label{lem:Jacob}
Suppose that $\bo$ is differentiable at $v$. 
Then, $\partial \bo (v) = \gamma P^{\pi_v}$, where $\pi_v$ is a greedy policy w.r.t.~$v$.
\end{Lem}
\begin{proof}
Let us define the matrix $P^a \in \R^{n\times n}$ with entries $P^a(s,s^+) = \PP\left(s^+ | s,a\right)$, for every control action $a \in \set{A}$.
Fix $s \in \set{S}$, and observe that  
\begin{align*}
     \partial_v \big( [\bo (v)](s) \big) 
     &= \partial_v \left( \min_{a\in\set{A}}\left\{ c(s,a) + \gamma\ \EE_{s^+}[v(s^+)\mid (s,a)] \right\}\right) 
     = \partial_v \left( \min_{a\in\set{A}}\left\{ c(s,a) + \gamma \cdot P^a(s,\cdot) \cdot v \right\}\right),
\end{align*}
where $ P^a(s,\cdot)$ is the $s$-th row of $P^a$. 
Then, using the envelope theorem~\cite{samuelson1948foundations}, we have
\begin{align*}
     \partial_v \big( [\bo (v)](s) \big) 
     &= \partial_v \bigg( c\big(s,\pi_v(s)\big) + \gamma \cdot P^{\pi_v}(s,\cdot) \cdot v \bigg) 
     = \gamma \cdot P^{\pi_v}(s,\cdot). 
\end{align*}
Hence, $\partial\bo (v) = \gamma P^{\pi_v}$, and the claim follows.
\end{proof}

Using the preceding result, we have 
$$\partial (\id - \bo)(v) = I - \gamma P(v),$$ 
where $P(v) \Let P^{\pi_v}$ is the state transition probability matrix of the greedy policy $\pi_{v_k}$ w.r.t.~$v_k$. 
Then, comparing the characterizations \eqref{eq:fp_v} and \eqref{eq:fp_gradient}, we can draw the following equivalence relations between deterministic optimization and model-based control:
\begin{align*}
   &  x \leftrightarrow v, \quad
    \left\{\begin{array}{rcl}
         \gr f & \leftrightarrow & \id - \bo \\
         \id - \gr f & \leftrightarrow & \bo
    \end{array}\right., \quad
    \left\{\begin{array}{rcl}
         \gr^2 f & \leftrightarrow & I - \partial\bo = I - \gamma P \\
         I - \gr^2 f & \leftrightarrow & \partial\bo = \gamma P 
    \end{array}\right. ,
\end{align*}
where $\id$ is the identity operator, $I$ is the identity matrix, and $P = P(v)$ is the transition probability matrix of the Markov chain under the greedy policy w.r.t.~$v$. 
Similarly, for stochastic optimization and model-free control, the characterizations \eqref{eq:fp_q} and \eqref{eq:expected_fp_gradient} point to 
the following equivalence relations:
\begin{align*}
    &\ \ \ \  (x ,  \xi) \leftrightarrow ( q ,  \varsigma^+ ),\quad
    \left\{\begin{array}{rcl}
         \gr \hf & \leftrightarrow & \id - \sbo\\
         \id - \gr \hf & \leftrightarrow & \sbo 
    \end{array}\right.,\quad
    \left\{\begin{array}{rcl}
         \gr^2 \hf & \leftrightarrow & I - \partial\sbo = I - \gamma \wh{P}\\
         I - \gr^2 \hf & \leftrightarrow & \partial\sbo = \gamma \wh{P} 
    \end{array}\right. ,
\end{align*}
where $\wh{P} = \wh{P}(q,\varsigma^+) \in [0,1]^{(\set{S}\times\set{A})\times (\set{S}\times\set{A})}$ is the \emph{synchronously sampled} transition probability matrix of the Markov chain under $\pi_q$ with elements\footnote{Once again, strictly speaking, $\wh{P}(q, \varsigma^+)$ is the empirical version of the \emph{state-action} transition probability matrix $\bar{P}(q) \in [0,1]^{(\set{S}\times\set{A})\times (\set{S}\times\set{A})}$ of the Markov chain under the greedy policy $\pi_q$ w.r.t.~$q$, with elements $[\bar{P}(q)]\big((s,a),(s',a')\big) = \PP(s' | s, a)$ if $a' = \pi_q(s')$ and $=0$ otherwise, for $(s,a),(s',a') \in \set{S}\times\set{A}$.} 
\begin{equation*}\label{eq:P_sampled}
     [\wh{P}(q,\varsigma^+)]\big((s,a),(s',a')\big) = \left\{\begin{array}{ll}
        1, & s' = \varsigma^+(s,a),\ a' = \pi_q\big(\varsigma^+(s,a)\big), \\ 
        0, & \text{otherwise},
     \end{array} \right. 
\end{equation*}
for each $\big((s,a),(s',a')\big) \in (\set{S}\times\set{A})\times(\set{S}\times\set{A})$.

\input{Tab/extended_algorithms_table.tex}

\section{Algorithms: Optimization vs. Control}\label{sec:equiv_alg} 
We now look at existing algorithms for optimization and control and their equivalence within the proposed framework. 
We start with the well-known baseline first-order and second-order algorithms and show how the application of the proposed transformations to well-established optimization algorithms, such as gradient descent and Newton method, leads to well-known control algorithms such as value iteration and policy iteration. 
We then cover a wide range of modifications to optimization algorithms, from momentum to spectral decomposition of the Hessian, and their counterpart control algorithms. 
We note that most of these equivalences have already been pointed out in the existing literature, however, mostly in a scattered way. 
An exception is~\cite{grand2021convex}, which studies the relation between \emph{deterministic} optimization algorithms and \emph{model-based} control algorithms.

To ease the exposition, we unify iterative algorithms by 
\begin{align*}
    y_{k+1} = y_k + d_k,\quad k = 0,1, \ldots,
\end{align*}
where $y_k = x_k$, $v_k$ or $q_k$ based on the context. 
In all considered settings, $d_k$ represents the update vector between iterations $k$ and $k+1$.
This form allows us to characterize algorithms in terms of $d_k$. 
Moreover, we use Greek letters $\alpha_k, \beta_k, \ldots$ as scalar coefficients (e.g., step-size, learning rate, etc.) in the update vectors. 
We note that the specific values of these parameters may differ between equivalent algorithms. 
A compact summary of this can be found in Table~\ref{tab:algorithms}. 

Some remarks are in order regarding the following equivalence relationships. 
First, the following list of algorithms is not meant to be an exhaustive review of the literature. Instead, our aim is to provide illustrative examples of the application of transformations of Table~\ref{tab:equivalence} in identifying equivalent algorithms for optimization and control. 
Second, the provided update rules for the cited algorithms, although covering the main characteristic of interest, are occasionally not exact. 
This is particularly the case in stochastic/model-free settings where other techniques, such as batch estimation or variance reduction, are also used in the original algorithms. 
Finally, the equivalence between optimization and control algorithms does not imply that these algorithms enjoy similar theoretical guarantees, particularly when it comes to convergence. 
We will discuss this issue in more detail in Section~\ref{sec:convergence}.   

\subsection{Baseline methods} 

We start our discussion with the vanilla first-order and second-order algorithms.

\textbf{Vanilla first-order methods.} 
The celebrated \emph{Gradient Descent (GD)}~\cite{lemarechal2012cauchy} method is characterized by 
$$d_k = -\alpha_k \gr f(x_k),$$ 
where $\alpha_k$ is a properly chosen step-size. 
Applying the transformations of Table~\ref{tab:equivalence} to GD, we derive the so-called \emph{Relaxed Value Iteration (Rel-VI)}~\cite{kushner1971accelerated,porteus1978accelerated,goyal2019first} with 
$$d_k = - \alpha_k \big(v_k - \bo (v_k)\big),$$ 
for model-based control. 
In particular, for the constant step-size $\alpha_k = 1$, we have the standard VI algorithm $v_{k+1} = \bo (v_k)$~\cite{bellman1957markovian}. 
Correspondingly, \emph{Stochastic GD (SGD)}~\cite{robbins1951stochastic} is characterized by 
$$d_k = -\alpha_k \gr \hf(x_k,\xi_k).$$ 
Under the transformations of Table~\ref{tab:equivalence}, SGD leads to the \emph{synchronous} \emph{Q-Learning (QL)} algorithm~\cite{watkins1992q,kearns1998finite} with\footnote{This is the so-called \emph{synchronous} update of the Q-function in \emph{all} state-action pairs in each iteration, corresponding to the \emph{parallel sampling model} introduced by~\cite{kearns1998finite}.} 
\begin{equation*}\label{eq:QL}
    d_k = - \alpha_k \big(q_k - \sbo (q_k,\varsigma^+_k)\big).
\end{equation*}

\textbf{Vanilla second-order methods.} 
In second-order algorithms, $d_k$ is specified by both the gradient and the Hessian oracles.
The \emph{damped} \emph{Newton Method (NM)} is one such algorithm with  
$$d_k = - \alpha_k [\gr^2f(x_k)]^{-1}\gr f(x_k).$$ 
The \emph{pure} Newton step with $\alpha_k = 1$ has a \emph{local} quadratic convergence if, in addition to $f$ being strongly convex, the Hessian is Lipschitz-continuous~\cite[Thm.~5.3]{bubeck2015convex}. 
Globally, however, the pure Newton method can lead to divergence. 
This is the reason behind introducing the step-size $\alpha_k < 1$ in the damped version which can be used to guarantee a global linear convergence.  
We can use the transformations of Table~\ref{tab:equivalence} in order to transform the Newton method into a model-based control algorithm with  
$$d_k = - \big(I - \gamma P(v_k)\big)^{-1}(v_k - \bo (v_k) ),$$ 
where $P(v_k)$ is the transition probability matrix of the Markov chain under the greedy policy w.r.t.~$v_k$. 
The derived model-based control algorithm then corresponds to the well-known \emph{Policy Iteration (PI)} algorithm~\cite{puterman1979convergence} with 
$v_{k+1} = \big(I - \gamma P(v_k)\big)^{-1} c^{\pi_{v_k}}$, 
where $c^{\pi_{v_k}}$ is the vector of stage costs corresponding to the greedy policy~$\pi_{v_k}$ w.r.t.~$v_k$. 
Indeed, the PI algorithm is equivalent to the semi-smooth Newton method with a local quadratic convergence rate~\cite{gargiani2022dynamic}. 

In the model-free case, the stochastic version of the Newton method~\cite{ruppert1985newton} has been a source of inspiration for developing second-order-type Q-learning algorithms. 
In particular, the \emph{Stochastic Newton-Raphson (SNR)}~\cite{ruppert1985newton} algorithm with
\begin{equation*}
    \left\{\begin{array}{l}
         D_k = (1-\beta_k) D_{k-1} + \beta_k \gr^2\hf(x_k,\xi_k), \\
         d_k = -\alpha_k D_k^{-1} \gr \hf(x_k,\xi_k),
    \end{array}\right.
\end{equation*}
was used for developing the \emph{Zap QL (ZQL)} algorithm~\cite{DevrajMeyn2017NIPS} with
\begin{equation*}
    \left\{\begin{array}{ll}
         D_k = (1-\beta_k) \ D_{k-1} + \beta_k \ e_{(s_k,a_k)} \ \big(e_{(s_k,a_k)} - \gamma \ e_{(s^+_k,\pi_k(s^+_k))}\big)\tr,\\[1ex]
         d_k = -\alpha_k \ D_k^{-1} \ \big( q_k(s_k,a_k) - [\sbo (q_k,\varsigma^+_k)](s_k,a_k) \big) \ e_{(s_k,a_k)},
    \end{array}\right.
\end{equation*}
where $\pi_k(s^+_k) \in \argmin_{a\in \set{A}} q_k(s^+_k,a)$ is a greedy action w.r.t.~$q_k$ evaluated at the sampled next state $\varsigma^+_k(s_k,a_k) = s^+_k\sim \PP(\cdot|s_k,a_k)$.
Note that the preceding algorithm involves updating one entry of the action-value function $q_k$ at each iteration $k$, corresponding to the state-action pair $(s_k,a_k)$ chosen at iteration~$k$ 
-- recall that $e_{(s,a)} \in \R^{\set{S}\times\set{A}}$ is the unit vector corresponding to the state-action pair~$(s,a)$.
Moreover, ZQL in~\cite{DevrajMeyn2017NIPS} implements an eligibility trace with a decaying factor which we omit (i.e., set to $0$) for simplicity of the representation.
The implementation of zap QL algorithm with \emph{synchronous} update of the Q-function in all state-action pairs in each iteration 
is then characterized by
\begin{equation}\label{eq:ZQL}
    \left\{\begin{array}{l}
        D_k = (1-\beta_k) D_{k-1} + \beta_k  \big(I - \gamma \wh{P}(q_k,\varsigma^+_k) \big), \\
        d_k = -\alpha_k D_k^{-1} \big( q_k - \sbo (q_k,\varsigma^+_k)\big),
    \end{array}\right.
\end{equation}
where $\wh{P}(q,\varsigma^+)$ is the synchronously sampled state-action transition probability matrix of the Markov chain under the greedy policy w.r.t.~$q$. 
Note that \eqref{eq:ZQL} is exactly the SNR algorithm under the transformations of Table~\ref{tab:equivalence}. 

\subsection{Accelerated methods}

We now focus on classic acceleration techniques such as momentum and anchoring in optimization. These techniques have recently attracted a lot of attention in developing control algorithms with improved convergence and complexity properties. 
We also look at a classic control engineering approach, namely, proportional-integral-derivative (PID) control, and its application in optimization and control algorithms.  

\textbf{Momentum methods.} 
In the so-called momentum-based algorithms, the update vector $d_k$ is specified by gradient oracles but also depends on $d_{k-1}$. 
One such algorithm is \emph{GD with Polyak Momentum (Mom-GD)}~\cite{polyak1964some}, a.k.a.~heavy ball method, characterized by 
$$d_k = -\alpha_k \gr f(x_k) + \beta_k d_{k-1}.$$  
Another well-known momentum-based algorithm is \emph{GD with Nesterov Acceleration (Acc-GD)}~\cite{nesterov1983method} with update vector
\begin{equation}\label{eq:Acc-GD}
    d_k = -\alpha_k \gr f(x_k+\beta_k d_{k-1}) + \beta_k d_{k-1}.
\end{equation}
With a proper choice of the step-sizes $\alpha_k$ and $\beta_k$, these schemes can be shown to accelerate the convergence rate, compared to the standard GD, for particular classes of objective functions~\cite{polyak1964some,nesterov2018lectures}. 
The corresponding model-based control algorithms, using the transformations of Table~\ref{tab:equivalence}, are \emph{Momentum VI (Mom-VI)}~\cite{goyal2019first} with 
\begin{equation}\label{eq:Mom-VI}
    d_k = - \alpha_k \big(v_k - \bo (v_k)\big) + \beta_k d_{k-1},
\end{equation}
and \emph{Accelerated VI (Acc-VI)}~\cite{goyal2019first} with 
$$d_k = - \alpha_k \big(v_k+\beta_k d_{k-1} -\bo (v_k+\beta_k d_{k-1})\big) + \beta_k d_{k-1}.$$ 
However, the convergence of the preceding accelerated schemes is in general not guaranteed. 
In~\cite{goyal2019first}, the authors address this issue by \emph{safeguarding}, i.e., combining the accelerated VI with the standard VI; see Section~\ref{sec:convergence} for more details. 

For accelerating SGD, a direct combination of Polyak momentum or Nesterov acceleration with SGD has been shown to lead to no better (and even worse) performance in terms of convergence rate~\cite{yang2016unified, kidambi2018insufficiency}. 
At least for almost sure convergence,~\cite{liu2022almost} reports the same rate of convergence for SGD with Polyak momentum and SGD with Nesterov acceleration as for standard SGD.
Nevertheless, modifications of momentum-based acceleration methods have led to a range of accelerated SGD algorithms with faster convergence rates with specific assumptions on the problem data~\cite{kidambi2018insufficiency,liu2018accelerating,allen2017katyusha}.  
The idea of using momentum for accelerating QL has also attracted some interest. 
In particular, applying the transformations of Table~\ref{tab:equivalence} on a generic \emph{Momentum SGD (Mom-SGD)}~\cite{yang2016unified} with 
\begin{equation*}
    \left\{\begin{array}{l}
         d'_{k-1} =  \gr \hf (x_{k},\xi_k) - \gr \hf (x_{k-1},\xi_k),\\
         d_k = -  \alpha_k \gr \hf (x_{k},\xi_k) - \beta_k d'_{k-1} + \delta_k d_{k-1},
    \end{array}\right.
\end{equation*}
and step-sizes $\alpha_k,\beta_k, \delta_k$, we obtain the \emph{Speedy QL (SQL)}~\cite{ghavamzadeh2011speedy}, \emph{Nesterov Stochastic Approximation (NeSA)}~\cite{devraj2019matrix}, and \emph{Momentum QL (Mom-QL)}~\cite{weng2020momentum} algorithms with 
\begin{equation}\label{eq:SQL}
    \left\{\begin{array}{l}
         d'_{k-1} =  \big(q_k -\sbo (q_{k},\varsigma^+_k)\big) - \big(q_{k-1} -\sbo (q_{k-1},\varsigma^+_k)\big) ,\\
         d_k = -  \alpha_k \big(q_k -\sbo (q_{k},\varsigma^+_k)\big) - \beta_k d'_{k-1} + \delta_k d_{k-1}. 
    \end{array}\right.
\end{equation}
The difference between these three algorithms is in the choice of the step-sizes $\alpha_k, \beta_k, \delta_k$. 
For example, by setting $\beta_k = \delta_k = (1 - 2 \alpha_k) = (k-1)/(k+1)$, we recover SQL. 

\textbf{Anchoring \& Halpern iteration.} 
Halpern iteration~\cite{halpern1967fixed} is an acceleration scheme originally developed for fixed-point iterations involving \emph{non-expansive} maps. 
In particular, it modifies the updates by gradually pulling the iterates toward a reference point, a.k.a.~\emph{anchor}. 
In the case of deterministic optimization, modifying the standard GD using Halpern iteration, we derive the \emph{Anchored GD (Anc-GD)}~\cite{tran2022connection}
\begin{equation*}
x_{k+1} = \beta_k x + (1 - \beta_k)\big( x_k -\delta_k \nabla f(x_k)\big),
\end{equation*}
where $x\in\R^n$ is the anchor, $\beta_k \in (0,1)$ is a properly chosen decaying step-size, and $\delta_k$ is the step-size corresponding to the standard GD update. 
The canonical choices for the anchor and the step-size are $x=x_0$ and $\beta_k=\tfrac{1}{k+2}$~\cite{lieder2021convergence}. 
Using our standard characterization, Anc-GD (with the anchor $x = x_0$) corresponds to the update vector
\[
d_k = \beta_k (x_0-x_k) - \alpha_k \nabla f(x_k),
\]
where $\alpha_k = \delta_k (1 - \beta_k)$. 
Under certain conditions on the objective function~$f$ and the step-sizes~$\alpha_k$ and $\beta_k$, Anc-GD improves the rate of convergence compared to standard GD. 
Indeed, Anc-GD can be shown to be equivalent to Nesterov acceleration with a correction term~\cite{tran2022connection}. 
However, when the operation corresponding to standard GD is a \emph{contraction}, the modification via anchoring in Anc-GD does \emph{not} improve the convergence rate~\cite{park2022exact}. 
The corresponding model-based control algorithm is \emph{Anchored VI (Anc-VI)}~\cite{lee2023accelerating} defined by the update vector
\begin{align*}
    \label{eq:anchorVI}
    d_k = \beta_k (v_0- v_k)  - (1- \beta_k) \big(v_k - \bo (v_k)\big),
\end{align*}
where $\alpha_k = (1- \beta_k)$.
A similar, undesired situation arises for Anc-VI algorithm. Since the Bellman operator is a contraction, there is no improvement in the convergence rate for $\gamma < 1$ and the iterates of Anc-VI converge to the optimal value function~$v\opt$ linearly with rate $\gamma$ (i.e., similar to standard VI). 
However, of importance is the polynomial (i.e., with rate~$\mathcal{O}(1/k)$) convergence of the Bellman residual~$\|v_k - \bo (v_k)\|_\infty$ to zero in the \emph{long-horizon} setting with $\gamma \ra 1$~\cite[Thm.~2]{lee2023accelerating}. 
In particular, even in the \emph{undiscounted} setting (i.e., $\gamma = 1$), the Bellman residual in Anc-VI converges to zero with rate~$\mathcal{O}(1/k)$, assuming the corresponding undiscounted Bellman operator has a fixed point~\cite[Thm.~3]{lee2023accelerating}.
As a result, in long-horizon and undiscounted settings, Anc-VI maintains robust convergence guarantees, while the geometric rate of standard VI deteriorates.

The corresponding \emph{Stochastic Halpern Iteration (SHI)}~\cite{cai2022stochastic,alacaoglu2025towards}, characterized by
\begin{equation*}\label{eq:shc}
    d_k = \beta_k (x_0-x_k) - \alpha_k \gr \hf(x_k,\xi_k),
\end{equation*}
with step-sizes $\alpha_k,\beta_k$, has been recently studied extensively in the stochastic optimization setting.
In particular, under certain conditions and in combination with other techniques such as variance reduction and restarting, anchoring leads to improved sample complexity~\cite{cai2022stochastic}. 
Moreover, anchoring has been shown to be instrumental in providing guarantees for the last iterate, as opposed to a weighted average of iterates, under relaxed bounded-variance assumptions~\cite{alacaoglu2025towards}. 
Applying the transformations of Table~\ref{tab:equivalence} on SHI, we derive the \emph{Halpern QL (HQL)}~\cite{bravo2024stochastic} and \emph{Stochastic Anchored Value Iteration for Discounted MDPs (SAVID)}~\cite{lee2025near} algorithms with 
\begin{equation*}\label{eq:hql}
    d_k = \beta_k (q_0-q_k) - (1- \beta_k) \big(q_k -\sbo (q_{k},\varsigma^+_k)\big),
\end{equation*}
where $\alpha_k = 1- \beta_k$. 
To be precise, HQL uses batch estimation with multiple samples for the sampled Bellman operation in this update rule to reduce the variance, while SAVID incorporates the recursive sampling technique~\cite{jin2024truncated} to reduce the sample complexity. 
Of importance is again that both HQL and SAVID algorithms (with minor modifications) can also be used in the \emph{averaged cost} setting in which the corresponding Bellman optimality operator is non-expansive but \emph{not} contractive. 

\textbf{PID control.} 
An alternative to standard acceleration methods is to use proportional-integral-derivative (PID) controllers for designing the update rule. 
This is a well-established technique in automatic control for set-point tracking. 
The basic idea is to treat the iterates~$x_k$ as the ``state'' of the system and define the ``error'' signal in terms of the gradient~$\gr f(x_k)$ and possibly the finite difference of the state, i.e., $d_k$.     
This idea has recently been explored in training deep networks in~\cite{an2018pid, chen2024accelerated}. 
In particular, the \emph{full batch} version of the \emph{PID Optimizer (PID-Opt)}~\cite{an2018pid} algorithm is a deterministic optimization algorithm characterized by 
\begin{equation*}
    \left\{\begin{array}{l}
        d^I_{k} = - \alpha_k \nabla f(x_k) + \beta_k d^I_{k-1},  \\[1ex]
        d^D_{k} =  \nabla f(x_k) - \nabla f(x_{k-1}), \\[1ex]
        d_k =  d^I_{k} +  d^D_{k},
        \end{array}\right.
\end{equation*}
where $\alpha_k$ is the step-size and $\beta_k$ is the decay coefficient in the integral term.  
The corresponding model-based control algorithm is \emph{PID accelerated VI (PID-VI)}~\cite{farahmand2021pid} with
\begin{equation*}
    \left\{\begin{array}{l}
        d^I_{k} = - \alpha_k \big(v_k - \bo (v_k)\big) + \beta_k d^I_{k-1},  \\[1ex]
        d^D_{k} = d_{k-1}, \\[1ex]
        d_k = - \kappa^P_k \big(v_k - \bo (v_k)\big) + \kappa^I_k d^I_{k}  + \kappa^D_k d^D_{k}.
        \end{array}\right.
\end{equation*}
Above, we see all three terms of the PID controller with the corresponding parameters $\kappa^P_k,\ \kappa^I_k,\ \kappa^D_k$. 
In particular, observe that, compared to PID-Opt, PID-VI uses the finite difference of the state (i.e., $d_{k-1}$) as opposed to the finite difference of the gradient (i.e., $d'_{k-1} = \nabla f(x_k) - \nabla f(x_{k-1})$) in its derivative term. 
Moreover, in~\cite{farahmand2021pid}, the authors provide an adaptive scheme for tuning the PID parameters. 
The close connection between the PID acceleration schemes above and the momentum-based algorithms can be readily seen in their update rules. 
Indeed, PID-Opt and PID-VI provided above can recover standard momentum-based algorithms for specific values of the PID parameters. 
For instance, Acc-GD~\eqref{eq:Acc-GD} can be derived by removing the derivative term in PID-Opt (i.e., setting $\kappa^D_k = 0$)~\cite{an2018pid}. 
Similarly, Mom-VI~\eqref{eq:Mom-VI} is an instance of PID-VI without the integral term  (i.e., with $\kappa^I_k = 0$)~\cite{farahmand2021pid}. 

For stochastic optimization, we can use the \emph{mini-batch} version of PID-Opt~\cite{an2018pid}, characterized by 
\begin{equation*}
    \left\{\begin{array}{l}
        d^I_{k} = - \alpha_k \nabla \hf(x_k) + \beta_k d^I_{k-1},  \\[1ex]
        d^D_{k} = \beta_k d^D_{k-1} + (1 - \beta_k) \big( \nabla \hf(x_k) - \nabla \hf(x_{k-1}) \big), \\[1ex]
        d_k = d^I_{k+1}  + \kappa_k^d d^D_{k},
        \end{array}\right.
\end{equation*}
where the stochastic gradient operation is evaluated using a mini-batch of samples.
In particular, observe that the derivative term also includes averaging to reduce the effect of noise in gradient evaluations. 
The corresponding model-free control algorithm is \emph{PID accelerated QL (PID-QL)}~\cite{bedaywi2024pid} with 
\begin{equation*}
    \left\{\begin{array}{l}
        d^I_{k} = - \alpha_k \big(q_k -\sbo (q_{k},\varsigma^+_k)\big) + \beta_k d^I_{k-1},  \\[1ex]
        q'_{k-1} = (1-\eta_k) q'_{k-2} + \eta_k q_{k-1}, \\[1ex]
        d^D_{k} = q_{k} - q'_{k-1}, \\[1ex]
        d_k = - \kappa^P_k \big(q_k -\sbo (q_{k},\varsigma^+_k)\big) + \kappa^I_k d^I_{k+1}  + \kappa^d_k d^D_{k},
        \end{array}\right.
\end{equation*}
In particular, observe that PID-QL also uses a stochastic approximation $q'_{k-1}$ of $q_{k-1}$ in forming the derivative term based on the finite difference of the state $q_k$. 
Moreover, a scheme similar to that used in PID-VI~\cite{farahmand2021pid} is used to tune the PID parameters in PID-QL adaptively. 
We also note that PID-QL can recover Mom-QL~\cite{weng2020momentum} and SQL~\cite{ghavamzadeh2011speedy} in~\eqref{eq:SQL} by setting the PID parameters accordingly~\cite{bedaywi2024pid}. 
Finally, while the empirical results reported for PID-Opt, PID-VI, and PID-QL show a significant improvement in the convergence rate compared to the corresponding standard methods, theoretical guarantees for the convergence of these algorithms are still lacking.   

\subsection{Quasi-Newton methods} 

Quasi-Newton methods (QNMs) are a class of methods that allow for a trade-off between computational complexity and (local) convergence rate in solving optimization problems. 
To do so, these methods use a Newton-type update vector
$$ d_k = - \alpha_k \wt{H}_k^{-1}\gr f(x_k),$$ 
where $\wt{H}_k$ (or $\wt{H}_k^{-1}$) is a computationally efficient approximation of the true Hessian $\gr^2f(x_k)$ (or its inverse) at iteration $k$. 
Different QNMs use different approximations of the Hessian or its inverse. 
Some ``classical'' examples of QNMs include Broyden approximation~\cite{broyden1965class}, Power symmetric Broyden approximation~\cite{powell1970new}, Davidon-Fletcher-Powell approximation~\cite{davidon1991variable,fletcher1963rapidly}. 
All of these examples approximate the Hessian (or its inverse) as the ``closest'' matrix to a prior matrix within an affine set defined by (multi-)secant constraints (and possibly other structural constraints such as symmetry, sparsity, etc.). 
Under certain conditions, these methods lead to a local \emph{superlinear} rate of convergence. 
Not surprisingly, different types of QNMs have been (directly or with some modifications) employed in solving model-based and model-free control problems. 
Some recent examples include the Sketched Newton VI~\cite{liu2024sketched} and Quasi-PI~\cite{kolarijani2023optimization} algorithms.  
In what follows, we focus on two illustrative examples, namely, Anderson mixing and spectral preconditioning. 

\textbf{Anderson mixing.} 
Anderson mixing, originally developed as an acceleration scheme for fixed-point iterations~\cite{anderson1965iterative}, provides an update rule in which the next iterate is constructed as a linear combination of a finite number of most recent iterates and their images. 
In particular, \emph{Anderson Accelerated GD (AA-GD)}~\cite{mai2020anderson}, with memory $m \ge 0$ and $m_k = \min\{k,m\}$, is characterized by the update vector
\begin{equation*}
    \left\{\begin{array}{l}
         X_k = [x_{k}, \ldots, x_{k-m_k}], \\[1ex]
         G_k = [\gr f(x_{k}), \ldots, \gr f(x_{k-m_k})],\\[1ex]
         w_k = (\e^\top (G_k^\top G_k)^{-1} \e)^{-1} (G_k^\top G_k)^{-1} \e,\\[1ex]
         d_k = - x_k + (X_k - \alpha_k G_k)w_k,
    \end{array}\right.
\end{equation*}
where $\e$ is the all-ones vector, $\alpha_k$ is a properly chosen step-size, and $G_k$ is assumed to be of full column rank (otherwise, we need to either use the pseudo-inverse or some form of regularization). 
Observe that for $m=0$, we recover the standard GD. 
Applying the transformations of Table~\ref{tab:equivalence}, we derive the \emph{Anderson Accelerated VI (AA-VI)}~\cite{geist2018anderson} with
\begin{equation}\label{eq:AA-VI}
    \left\{\begin{array}{l}
         V_k = [v_{k}, \ldots, v_{k-m_k}], \\[1ex]
         G_k = [v_{k} - \bo(v_k), \ldots, v_{k-m_k} - \bo(v_{k-m_k})],\\[1ex]
         w_k = (\e^\top (G_k^\top G_k)^{-1} \e)^{-1} (G_k^\top G_k)^{-1} \e,\\[1ex]
         d_k = - v_k + (V_k - G_k)w_k,
    \end{array}\right.
\end{equation}
where $\alpha_k = 1$ and $G_k$ is assumed to be of full column rank. 
The preceding update rules correspond to the so-called \emph{type-II} Anderson mixing. 
To be precise, the update vector of AA-GD/VI admits a QNM characterization of the form $d_k = - D_k g(y_{k})$, where $D_k$ is an approximation of the \emph{inverse} of the Hessian, i.e.,  $H_k^{-1}$~\cite{fang2009two}. 
In~\cite{zhang2020globally}, the authors consider the \emph{type-I} Anderson mixing scheme in combination with GD/VI in which the update vector is of the form $d_k = - D_k^{-1} g(y_{k})$ with $D_k$ being an approximation of the Hessian~$H_k$. 
We note that the combination of Anderson mixing (type-I or type-II) with VI also lacks convergence guarantees. In~\cite{zhang2020globally}, this issue is again addressed via \emph{safeguarding}. 

Anderson mixing has also been used in the stochastic optimization setting. 
In particular, in~\cite{NEURIPS2021_c203e4a1}, the authors use (type-II) Anderson mixing to develop the \emph{Stochastic Anderson Mixing (SAM)} algorithm with update vector
\begin{equation*}
    \left\{\begin{array}{l}
         D^{x}_k = [d^{x}_{k-m_k}, \ldots, d^{x}_{k-1}] \quad \text{with} \quad d^{x}_{i} = x_{i+1} - x_{i}, \\[1ex]
         D^{g}_k = [d^{g}_{k-m_k}, \ldots, d^{g}_{k-1}] \quad \text{with} \quad d^{g}_{i} = \gr \hf (x_{i+1},\xi_{i+1}) - \gr \hf (x_{i},\xi_i),\\[1ex]
         D_k = \beta_k I - \alpha_k(D^{x}_k + \beta_k D^{g}_k)\big((D^{g}_k)^\top D^{g}_k + R_k\big)^{-1}(D^{g}_k)^\top,\\[1ex]
         d_k = - D_k \gr \hf (x_{k},\xi_{k}),
    \end{array}\right.
\end{equation*}
where $\alpha_k, \beta_k$ are properly chosen step-sizes and $R_k = \delta_k (D^{x}_k)^\top D^{x}_k$, with an adaptive coefficient $\delta_k$, acts as a \emph{regularizer}. 
To be precise, SAM uses a batch of samples for estimation of the stochastic gradient operator. 
Note that the formulation above aligns with the QNM characterization of Anderson mixing. 
The corresponding model-free control algorithm is \emph{QL with Stable Anderson Acceleration (SAA-QL)}~\cite{sun2021damped} with
\begin{equation}\label{eq:SAA-QL}
    \left\{\begin{array}{l}
         D^{q}_k = [d^{q}_{k-m_k}, \ldots, d^{q}_{k-1}] \quad \text{with} \quad d^{q}_{i} = q_{i+1} - q_{i}, \\[1ex]
         D^{t}_k = [d^{t}_{k-m_k}, \ldots, d^{t}_{k-1}] \quad \text{with} \quad d^{t}_{i} = \sbo (q_{i+1},\varsigma^+_{i+1}) - \sbo (q_{i},\varsigma^+_{i}),\\[1ex]
         D^{g}_k = D^{q}_k - D^{t}_k, \\[1ex]
         D_k = \beta_k I - (D^{q}_k + \beta_k D^{g}_k)\big((D^{g}_k)^\top D^{g}_k + R_k\big)^{-1}(D^{g}_k)^\top,\\[1ex]
         d_k = - D_k \big(q_k -\sbo (q_{k},\varsigma^+_k)\big),
    \end{array}\right.
\end{equation}
where $\alpha_k = 1$ and $R_k = \delta_k (\|D^{q}_k\|_F^2 + \|D^{g}_k\|_F^2) I$ as the regularizer. 
To be precise, SAA-QL in~\cite{sun2021damped} utilizes a \emph{smoothed} version of the Bellman optimality operator based on the \emph{mellow-max} operator; see Section~\ref{subsec:model_free_SA} for more details. 
Moreover, in~\cite{sun2021damped}, the proposed SAA is integrated within the Dueling Deep Q-Network (Duel-DQN)~\cite{wang2016dueling} algorithm for \emph{continuous-state} model-free control. 
In particular, the SAA-Duel-DQN algorithm stores the last $m_k$ target networks (i.e., their parameters) and uses \eqref{eq:SAA-QL} for computing the target values in each iteration. 
Finally, let us note that several works have incorporated the original form of the Anderson mixing similar to~\eqref{eq:AA-VI} for model-free control; see, e.g.,~\cite{shi2019regularized,zuo2022offline}.

\textbf{Spectral preconditioning}. 
This class of QNMs involves another computationally efficient approximation of the Hessian. 
To be precise, they use an approximation of the form 
\[
\wt{H}_k = \alpha_k I + \beta_k U_kW_k\tr,
\]
where $U_k,W_k \in \R^{\ell\times r}$ are tall matrices with $r \ll \ell$. 
Such approximations readily allow for application of the Woodbury matrix identity~\cite{hager1989updating} for computing $\wt{H}_k^{-1}$ with a low computational cost. 
In spectral preconditioning techniques, the matrices $U_k$ and $W_k$ are particularly formed based on the spectral decomposition of the true Hessian. 
For instance, \emph{GD with Spectral Preconditioning (SP-GD)}~\cite{doikov2024spectral} is characterized by
\begin{equation*}
    \left\{\begin{array}{l}
         D_k = \alpha_k I + \sum_{i=1}^r \lambda_{k,i} u_{k,i} u_{k,i}\tr ,\\[1ex]
         d_k = - D_k^{-1} \gr f (x_k),
    \end{array}\right.
\end{equation*}
with $\alpha_k$ being a properly chosen regularizer parameter, $(\lambda_{k,i})_{i=1}^r$ the largest eigenvalues of the Hessian $\gr^2f(x_k)$ and $(u_{k,i})_{i=1}^r$ the corresponding eigenvectors. 
For model-based control, we can mention the recently developed \emph{Rank-One Modified VI (R1-VI)}~\cite{kolarijani2025rank} with 
\begin{equation*}
    \left\{\begin{array}{l}
         D_k = I - \gamma \e w_k\tr ,\\[1ex]
         d_k = - D_k^{-1} \big(v_k - \bo (v_k)\big),
    \end{array}\right.
\end{equation*}
where $w_k$ is the stationary distribution estimate of the Markov chain induced by the greedy policy w.r.t.~$v_k$. 
The preceding update rule particularly exploits the structure $H(v) = I-\gamma P(v)$ of the ``Hessian'' in the control problem and employs the spectral decomposition of the matrix $P(v)$. 
In this regard, we note that $w_k$ and the all-one vector~$\e$ are, respectively, the left and right eigenvectors of $P(v_k)$ corresponding to the largest eigenvalue, i.e., $1$. 

A similar idea has also been used in stochastic optimization. 
In particular, \emph{Newton-Sampling Method via Rank Thresholding (NewSamp)}~\cite{erdogdu2015convergence} is characterized by 
\begin{equation*}
    \left\{\begin{array}{l}
         D_k = \wh{\lambda}_{k,r+1} I + \sum_{i=1}^r (\wh{\lambda}_{k,i}-\wh{\lambda}_{k,r+1}) \wh{u}_{k,i} \wh{u}_{k,i}\tr ,\\[1ex]
         d_k = - D_k^{-1} \gr \hf (x_{k},\xi_{k}),
    \end{array}\right.
\end{equation*}
with eigenvalues $(\wh{\lambda}_{k,i})_{i=1}^{r+1}$ and the corresponding eigenvectors $(\wh{u}_{k,i})_{i=1}^r$ computed using eigenvalue decomposition of the \emph{batch-sampled} Hessian, that is, $\frac{1}{B} \sum_{j=1}^{B} \gr^2 \hf (x_k, \xi_{k,j})$. 
For model-free control, the same idea is used in  \emph{Rank-One Modified QL (R1-QL)}~\cite{kolarijani2025rank} with 
\begin{equation*}
    \left\{\begin{array}{l}
         D_k = I - \gamma \e \wh{w}_k\tr ,\\[1ex]
         d_k = - \alpha_k D_k^{-1} \big(q_k -\sbo (q_{k},\varsigma^+_k)\big),
    \end{array}\right.
\end{equation*}
where $\wh{w}_k$ is now the state-action stationary distribution of the Markov chain induced by the greedy policy w.r.t.~$q_k$ and estimated using the sample transitions. 
Once again, the preceding update rule exploits the structure of the ``Hessian'' in the control problem. 

We note that the spectral preconditioning technique is closely related to the \emph{matrix splitting} method for solving linear equations. 
In particular, Deflated Dynamics VI (DDVI) and Temporal Difference (DDTD) algorithms~\cite{lee2024deflated} also use a spectral decomposition of the transition probability matrix to form a low-rank approximation of it. They then use this low-rank approximation in combination with matrix splitting for the \emph{policy evaluation} problem. 
In effect, the corresponding update rules follow the same structure as R1-VI/QL, except with rank-$r$ approximation of the transition probability matrix where $r>1$ (as opposed to $r=1$). 

\section{Convergence of Control Algorithms}\label{sec:convergence}

One of the promises of building a common formulation for control and optimization algorithms is to borrow existing theoretical results from one field to another.
However, as several other studies have pointed out, the similarity in the update rule does not carry over to the analysis. 
One of the main reasons for this discrepancy is the \emph{non-differentiability} of the Bellman optimality operator~$\bo$. 
Another complicating factor is that the operator~$\bo$ is a contraction in $\infty$-norm, which is also \emph{non-differentiable}. 
This forces the convergence analysis of the control algorithms to rely more on fixed-point theory and contraction mappings than on gradient-based smooth analysis. 

In what follows, we examine general tools for the convergence analysis of control algorithms, as well as some simple yet effective remedies for ensuring the convergence of novel control algorithms. 

\subsection{Model-based algorithms}
\label{subsec:model_based}

For \emph{finite} state-action MDPs, the non-differentiability of the operator~$\bo$ is ``manageable.'' 
Indeed, for such MDPs, $\bo$ is piece-wise affine and we have  
\begin{align*}
\bo (v) = \min_{\pi\in \set{S}^{\set{A}}} c^{\pi} + \gamma P^{\pi}v, \quad \forall v \in \R^{\set{S}},
\end{align*}
where the minimization is over the finite set of all deterministic control policies, each corresponding to a linear function $c^{\pi}+\gamma P^{\pi} v$.
Equivalently, the finite family of policies admits an affine switching system \cite{liberzon2003switching}, a connection leveraged by~\cite{lee2020unified} to establish convergence guarantees for asynchronous Q-learning.
Moreover, by Rademacher's Theorem~\cite[Thm.~10.8]{villani2008optimal}, $\bo$ is differentiable almost everywhere. 
This property allows one to use tools from semismooth analysis for deriving local convergence results. 
In particular, in~\cite{gargiani2022dynamic}, the authors use the fact that the PI algorithm is an instance of the semi-smooth Newton method to show that it has a local quadratic convergence rate. Similarly, in~\cite{kolarijani2023optimization}, the semismoothness of $\bo$ is used for deriving a local superlinear convergence rate for one of the proposed quasi-PI algorithms.   

A more general approach, however, is to completely get rid of the non-differentiability by \emph{smoothing} the Bellman optimality operator. 
Not surprisingly, a similar technique has been extensively used to develop algorithms for solving non-smooth optimization problems; see, e.g., the seminal work~\cite{nesterov2005smooth}. 
The smoothing of $\bo$ is achieved by approximating the minimization operation by a differentiable function. 
As an example, by using the \emph{log-sum-exp} (a.k.a.~\emph{soft-max}) function~\cite{rust1994structural}, we can define
\begin{equation*}
\begin{array}{c}
[\bar{\bo}_{\beta} (q)](s,a) \Let c(s,a) + \gamma\ \EE_{s^+} \big[ - \frac{1}{\beta} \log \big( \ssum_{a^+ \in\set{A}} \text{e}^{-\beta q(s^+,a^+)}  \big) \;\big|\; (s,a)\big],\quad \forall (s,a)\in\set{S}\times\set{A},
\end{array}
\end{equation*}
where $\beta > 0$ is the so-called temperature parameter. 
The operator $\bar{\bo}_{\beta}$ is an over-approximation of the true Bellman optimality operator~$\bar{\bo}=\EE_{\varsigma^+} [\sbo (\cdot,\varsigma^+)]$ and, in particular, $\bar{\bo}_{\beta} \ra \bar{\bo}$ as $\beta \ra \infty$.
Accordingly, its fixed-point~$q^{\beta}$ is not the optimal Q-function, while the gap $\|q^{\beta} - q\opt\|_\infty$ can be controlled via the temperature parameter~$\beta$. 
However, such smoothing allows for a direct application of the second-order scheme for finding the fixed-point of $\bar{\bo}_{\beta}$~\cite{rust1994structural}. 
This idea has been recently used to propose the generalized second-order VI with a quadratic convergence rate~\cite{kamanchi2021generalized}. 

Besides the general tools discussed above, there are simple yet efficient ``safeguarding'' techniques for ensuring the convergence of novel model-based algorithms that lack convergence guarantees. 
These techniques specifically exploit the fact that the Bellman optimality operator~$\bo$ is a $\gamma$-contraction. 
To be precise, let us consider a generic model-based update rule 
\begin{align*}
v_{k+1} = v_k + d_k, 
\end{align*}
characterized by the novel update vector~$d_k$. 
One simple technique for ensuring the convergence of the preceding update rule is to safeguard it against VI~\cite{goyal2019first}. 
For some $\gamma' \in [\gamma,1)$, we consider the \emph{safeguarded update rule}
\begin{equation}\label{eq:safegaurd}
\begin{array}{l}
    v_{k+1} = v_k + d_k; \\[1.5ex]
    \text{{\bf if} $\|v_{k+1} - \bo (v_{k+1})\|_{\infty} > (\gamma')^k\ \|v_{0} - \bo (v_{0})\|_{\infty}$} \\[1.5ex]
    \hspace{1cm} v_{k+1} = \bo(v_k); \\[1.5ex]
    \text{{\bf endif}}
\end{array}
\end{equation}
The preceding update rule has the following convergence guarantee (for completeness, we provide the proof in Appendix~\ref{app:proofsafeguardVI}; see also~\cite[Thm.~4]{goyal2019first}):

\begin{Thm}[Convergence via safeguarding against VI]\label{thm:safeguardVI}
Let $\gamma' \in [\gamma,1)$. The safeguarded update rule~\eqref{eq:safegaurd}, with arbitrary update vector $d_k$, converges linearly to the optimal value function $v\opt = \bo(v\opt)$ with rate $\gamma'$. 
\end{Thm}

Alternatively, we can modify the update rule and use backtracking to ensure the convergence~\cite{kolarijani2023optimization}. 
Specifically, we modify the update rule as 
\begin{align*}
v_{k+1} = \bo(v_k) + \alpha_k \big(v_k - \bo(v_k) + \beta_k d_k\big),
\end{align*}
by introducing the coefficients $\alpha_k$ and $\beta_k$. Observe that we have the original update rule by setting $\alpha_k = \beta_k = 1$. 
The coefficient $\alpha_k$ is determined via backtracking to ensure that the Bellman error $\|v_{k+1} - \bo (v_{k+1})\|_{\infty}$ is reduced in each iteration linearly with some rate $\gamma' \in (\gamma,1)$. 
The coefficient $\beta_k$, on the other hand, is introduced to ensure that the novel update vectors are bounded and hence the backtracking terminates in finitely many steps. 
To be precise, we consider the \emph{update rule with backtracking} 
\begin{equation}\label{eq:backtrack}
\begin{array}{ll}
    \text{(0)} & \alpha_k = 1; \\[1.5ex]
    \text{(1)} & v_{k+1} = \bo(v_k) + \alpha_k \big(v_k - \bo(v_k) + \beta_k d_k\big); \\[1.5ex]
    \text{(2)} & \text{{\bf if} $\|v_{k+1} - \bo (v_{k+1})\|_{\infty} > \gamma'\ \|v_{k} - \bo (v_{k})\|_{\infty}$} \\[1.5ex]
    & \hspace{1cm} \text{$\alpha_{k} \gets \lambda \alpha_{k}$ and go to step~(1) }; \\[1.5ex]
    & \text{{\bf endif}}
\end{array}
\end{equation}
The preceding update rule has the following convergence properties (for completeness, we provide the proof in Appendix~\ref{app:proofbacktrackVI}; see also~\cite[Sec.~III.B.2]{kolarijani2023optimization}):

\begin{Thm}[Convergence via backtracking against VI]\label{thm:backtrackVI} 
Let $\gamma' \in (\gamma,1)$, $\lambda \in (0,1)$, and 
\begin{align*}
\beta_k = \frac{\min\{\norm{v_{k} - \bo (v_{k})}_{\infty},\ \norm{d_k}_\infty\}}{\norm{d_k}_\infty}, \quad \forall k \geq 0. 
\end{align*}
Then, the backtracking scheme in~\eqref{eq:backtrack} is guaranteed to terminate in $\ord\big(\ln(\gamma'-\gamma)/\ln(\lambda)\big)$ steps and the corresponding update rule, with arbitrary update vector $d_k$, converges linearly to the optimal value function $v\opt = \bo(v\opt)$ with rate $\gamma'$. 
\end{Thm}

Observe that the safeguarding and backtracking methods described above lead to a convergence behavior (i.e., linear with rate $\gamma'$) and per-iteration computational complexity (i.e., $\ord(n^2m)$ discarding logarithmic factors) similar to the standard VI algorithm. 
This is not surprising considering the fact that these methods, similar to VI, exploit the contractivity of the Bellman optimality operator~$\bo$. 
Nevertheless, these methods enable the safe integration of novel, optimization-inspired update rules (encapsulated in $d_k$) into a provably convergent control algorithm~\cite{goyal2019first, kolarijani2023optimization}.

\subsection{Model-free algorithms}
\label{subsec:model_free_SA}

Similar to a model-based case, one can benefit from smoothing approaches in model-free scenarios. 
In simple words, smoothing out the kinks in the Bellman optimality operator or the optimal policies allows us to borrow well-developed tools for stochastic optimization of smooth functions. 
The recently proposed \emph{mellow-max} operator~\cite{AsadiLittman2017AltSoftmax} is another excellent example of such smoothing that can be used as a differentiable approximation of the minimization within the Bellman optimality operator. 
Besides being a smooth operator, the mellow-max operator is, more importantly, non-expansive in the $\infty$-norm.
This property then streamlines the usage of standard convergence techniques seamlessly.
The corresponding smoothed Bellman optimality operator using the mellow-max is defined as
\begin{equation*}
\begin{array}{c}
[\bo_{\omega} (q)](s,a) \Let c(s,a) + \gamma\ \EE_{s^+} \big[ - \frac{1}{\omega} \log \big( \frac{1}{|\set{A}|} \sum_{a^+ \in\set{A}} \text{e}^{-\omega q(s^+,a^+) }  \big) \;\big|\; (s,a)\big],\quad \forall (s,a)\in\set{S}\times\set{A},
\end{array}
\end{equation*}
with $\omega > 0$ being the temperature parameter that can be used to control the gap $\|q^\omega-q\opt\|_\infty$ between the fixed-point~$q^\omega$ of $\bo_\omega$ and the true optimal Q-function $q\opt = \bar{\bo} (q\opt)$. 
The smoothed operator $\bo_\omega$ is an under-approximation of the true Bellman optimality operator~$\bar{\bo}$ and, in particular, we have $\bo_\omega \ra \bar{\bo}$ as $\omega \ra \infty$.  
The operator~$\bo_\omega$ has been recently used in combination with Anderson mixing to improve the convergence behavior in policy iteration~\cite{sun2021damped}.

In the model-free case, algorithmic updates are inherently noisy, being driven by samples of the dynamics and costs.
This stochasticity naturally places these algorithms within the scope of \emph{stochastic approximation} (\emph{SA}).
A standard SA iteration is
\begin{equation}\label{eq:SA_discrete}
    x_{k+1} = x_k + \alpha_k \big( h(x_k) + M_{k+1} \big),
\end{equation}
where $\alpha_k$ is the step-size (or the so-called \emph{learning rate}), $h:\R^\ell\!\to\!\R^\ell$ is the \emph{mean field} (expected update), and $M_{k}$ is a martingale-difference noise.
A common approach to studying the limiting behavior of an SA method is the \emph{ordinary differential equation} (\emph{ODE}) method.
The ODE method (with a dynamical-systems viewpoint), originating with Ljung~\cite{Ljung1977}, analyzes \eqref{eq:SA_discrete} by comparing it with the flow of
\begin{equation*}\label{eq:mean_field_ODE}
    \dot{x}(t) = h\big(x(t)\big),
\end{equation*}
and drawing conclusions about the iterates from the asymptotics of this ODE.
The interested reader is referred to \cite{Chen2002StochasticApproximation,KushnerYin2003,Borkar2008} for comprehensive treatments.
Invariance principles from nonlinear systems (e.g., LaSalle's invariance principle; see \cite[Ch.~4]{Khalil2002}) then identify invariant/limit sets of the ODE toward which trajectories converge~\cite{Benaim1996,Benaim1999}.
Under standard SA assumptions, the SA iterates track these sets.
To guarantee almost-sure boundedness (or Lyapunov stability of the iterates), one imposes verifiable drift/Lyapunov conditions. 
See \cite{BorkarMeyn2000} for the ODE method with drift at infinity and \cite{AndrieuMoulinesPriouret2005} for global Lyapunov constructions. 
It is worth mentioning that an important RL-specific deviation from otherwise standard SA conditions is \emph{ergodicity of the state process} for policy evaluation and \emph{sufficient exploration} for control.

Over the years, the ODE framework, in particular the framework introduced by~\cite{BorkarMeyn2000}, has become a standard tool for establishing stability and convergence of a broad class of model-free control/RL algorithms.
Considering the \emph{average cost problems}, in~\cite{AbounadiBertsekasBorkar2001}, the authors provided for the first time a complete analysis of two Q-learning algorithms by coupling the ODE style analysis with Kushner-Clark argumentation style~\cite{KushnerClark1978}. 
Another major subclass of RL methods is the well-known \emph{actor-critic} methods, which can be characterized as \emph{coupled}, \emph{two-timescale} stochastic iterations: the critic evolves on a \emph{fast} scale to track the value of the current policy, and the actor evolves on a \emph{slow} scale using the critic's signal to adjust the policy. 
The ODE viewpoint makes this decomposition explicit, which, under the usual SA conditions together with the RL ones, yields clean convergence guarantees to a stationary policy~\cite{BhatnagarEtAl2009NAC, PerkinsLeslie2012}. 
Another example is the use of the ODE method to show the convergence of \emph{gradient temporal difference} (\emph{TD}) \emph{approaches} under the general setup of off-policy learning and linear function approximation~\cite{SuttonEtAl2009ICML}. 
This work has then been extended to any smooth value function approximation in~\cite{MaeiEtAl2009NIPS} and to a control setting (policy improvement) under a stationary behavioral policy in~\cite{MaeiEtAl2010ICML}. 
Recently, ODE-style arguments have been used not only to prove convergence but also to \emph{design} RL algorithms. 
In~\cite{DevrajMeyn2017NIPS}, the authors introduce a \emph{matrix preconditioning} for the TD direction, improving the transient behavior of the proposed algorithm. This matrix gain is learned on a faster timescale and tracks the inverse Jacobian of the mean field.
This design extends to nonlinear function approximation under standard regularity assumptions~\cite{ChenEtAl2019arXiv}.

As an example of application of results developed using ODE methods for SAs, namely~\cite{PerkinsLeslie2012}, we finish this section by providing a simple and efficient ``safeguarding'' technique for ensuring the convergence of novel model-free algorithms that lack convergence guarantees. 
To this end, let us consider a generic update rule
\begin{align*}
q_{k+1} = q_k + \alpha_k b_k, 
\end{align*}
characterized by the novel update vector~$d_k = \alpha_k b_k$. 
Here, $\alpha_k$ is the learning rate satisfying the standard Robbins–Monro condition (i.e., $\sum_k \alpha_k = \infty$ and $\sum_k \alpha_k^2 < \infty$) and $b_k$ is a vector depending on the sampled Bellman error, a.k.a.~temporal difference (i.e., $\sbo (q_k,\varsigma^+_k)-q_k$). 
To address the convergence issue, one possibility is to consider a convex combination of the standard QL algorithm update rule~\eqref{eq:QL} with this novel update rule 
with a diminishing effect~\cite{kolarijani2023optimization}. 
To be precise, we modify the update rule as follows
\begin{align*}
q_{k+1} &= q_k + \alpha_k \left( (1-\beta_k) \big(\sbo (q_k,\varsigma^+_k)-q_k\big) + \beta_k  b_k \right)\\
&= q_k + \alpha_k \left(\sbo (q_k,\varsigma^+_k)- q_k + \beta_k \big(b_k + q_k-\sbo (q_k,\varsigma^+_k)\big) \right),
\end{align*}
where $\beta_k$ is a diminishing coefficient (i.e., $\beta_k \ra 0$). 
Finally, to ensure convergence, we need to bound the effect of the extra term $p_k = b_k + q_k -\sbo (q_k,\varsigma^+_k)$.
For that, given some $\rho > 0$, we employ the operator~$B_\rho: \R^{\set{S}\times\set{A}}\ra\R^{\set{S}\times\set{A}}$ defined by
\begin{equation*}
  B_\rho(p) \Let \frac{\min\{\rho, \norm{p}_\infty\}}{ \norm{p}_\infty} p.
\end{equation*} 
This leads to the \emph{safeguarded update rule}
\begin{equation}\label{eq:combined}
\begin{array}{l}
     p_k = b_k - q_k + \sbo (q_k,\varsigma^+_k),  \\[1ex] 
     q_{k+1} = q_k + \alpha_k \big(\sbo (q_k,\varsigma^+_k)-q_k + \beta_k B_\rho(p_k)\big), 
\end{array}
\end{equation}
with the following convergence property (for completeness, we provide the proof in Appendix~\ref{app:proofsafeguardQL}; see also~\cite[Sec.~III.C]{kolarijani2023optimization}): 

\begin{Thm}[Convergence via safeguarding against QL]\label{thm:safeguardQL}
Consider the safeguarded update rule~\eqref{eq:combined}.  
Let $\rho > 0$ be given and assume the learning rates $\alpha_k$ and $\beta_k$ are such that $\sum_k \alpha_k = \infty$, $\sum_k \alpha_k^2 < \infty$, and $\beta_k \ra 0$. 
Then, the iterates~$q_k$ converge to $q\opt$ almost surely. 
\end{Thm}

We note that the safeguarding technique described above leads to the same asymptotic convergence behavior as QL due to the diminishing effect of the extra term $p_k$ (which includes the novel direction~$b_k$). 
However, it does not involve any extra sample complexity since the computation of the vector $b_k$ already requires us to have access to the sampled Bellman error~$\sbo (q_k,\varsigma^+_k)-q_k$. 

\appendix

\section{Technical Proofs}\label{app:proofs}

\subsection{Proof of Theorem~\ref{thm:safeguardVI}}\label{app:proofsafeguardVI}

We will use induction to show that 
\begin{equation}\label{eq:proof_induction}
\norm{v_{k} - \bo(v_k)}_{\infty} \leq (\gamma')^{k} \norm{v_{0} - \bo(v_{0})}_{\infty}, \quad \forall k \geq 0.
\end{equation}
The base case $k=0$ holds trivially.
Let us now assume the inequality \eqref{eq:proof_induction} holds for some $k \geq 0$ and recall that $\bo$ is a $\gamma$-contraction in $\infty$-norm. 
Then, the safeguarding against standard VI as in \eqref{eq:safegaurd} implies that 
\begin{align*}
        \norm{v_{k+1} - \bo(v_{k+1})}_{\infty}  &\leq \max \{ (\gamma')^{k+1} \norm{v_{0} - \bo(v_{0})}_{\infty},\; \gamma \norm{v_{k} - \bo(v_{k})}_{\infty} \} \\
        &\leq \max \{ (\gamma')^{k+1} \norm{v_{0} - \bo(v_{0})}_{\infty},\; \gamma' \norm{v_{k} - \bo(v_{k})}_{\infty} \}\\
        &\leq  (\gamma')^{k+1} \norm{v_{0} - \bo(v_{0})}_{\infty},
\end{align*}
where we used $\gamma'\geq \gamma$ for the second inequality and the induction hypothesis for the last inequality. This completes the proof.

\subsection{Proof of Theorem~\ref{thm:backtrackVI}}\label{app:proofbacktrackVI}

Note that if $\|v_{k} - \bo (v_{k})\|_{\infty}=0$, we clearly have $v_k=v\opt$. 
Hence, we can assume $v_k\neq v\opt$. 
Then, using~\cite[Lem.~3.3]{kolarijani2023optimization}, we have that $\|v_{k+1} - \bo (v_{k+1})\|_{\infty} \leq \gamma'\ \|v_{k} - \bo (v_{k})\|_{\infty}$ if
\begin{align*}
\norm{v_{k+1} - \bo(v_k)}_\infty = \alpha_k \norm{v_k - \bo(v_k) + \beta_k d_k}_\infty \leq \frac{\gamma'-\gamma}{1+\gamma} \norm{v_{k} - \bo(v_k)}_\infty
\end{align*}
That is, backtracking terminates after $\ell \geq 0$ steps if 
\begin{align*}
\alpha_k^{-1} = \lambda^{-\ell} \geq \frac{1+\gamma}{\gamma'-\gamma} \cdot \frac{\norm{v_k - \bo(v_k) + \beta_k d_k}_\infty}{\norm{v_k - \bo(v_k)}_\infty}.
\end{align*}
Now, observe that $1+\gamma \leq 2$ and the choice of $\beta_k$ implies that 
\begin{align*}
\norm{v_k - \bo(v_k) + \beta_k d_k}_\infty \leq \norm{v_k - \bo(v_k) }_\infty + \beta_k \norm{ d_k}_\infty \leq 2 \norm{v_k - \bo(v_k) }_\infty.
\end{align*}
This, in turn, implies that the backtracking terminates after $\ell \leq L$ steps, where 
\begin{align*}
\lambda^{-L} = \frac{4}{\gamma'- \gamma} \iff L = \ord\big(\ln(\gamma'-\gamma)/\ln(\lambda)\big).
\end{align*}
The linear convergence of the iterates $v_k$ of update rule with backtracking to $v\opt$ with rate $\gamma'$ then immediately follows from the condition $\|v_{k+1} - \bo (v_{k+1})\|_{\infty} \leq \gamma'\ \|v_{k} - \bo (v_{k})\|_{\infty}$, which implies
\begin{align*}
\norm{v_{k} - \bo(v_k)}_{\infty} \leq (\gamma')^{k} \norm{v_{0} - \bo(v_{0})}_{\infty}, \quad \forall k \geq 0.
\end{align*}

\subsection{Proof of Theorem~\ref{thm:safeguardQL}}\label{app:proofsafeguardQL}

Define 
\begin{align*}
    h(q) \Let \bar{\bo}(q) - q,\
    M_k \Let \sbo (q_k,\varsigma^+_k) - \bar{\bo}(q_k),\
    r_k \Let \beta_k B_\rho(p_k),
\end{align*}
where $\bar{\bo}$ is the corresponding true Bellman optimality operator for Q-functions, that is, $\bar{\bo}(q) = \EE_{\varsigma^+} [\sbo (q,\varsigma^+)]$ for each $q$. 
Then, the convergence of $q_k$ to the unique equilibrium $q\opt$ of the map $h$ follows from \cite[Cor.~3.2]{PerkinsLeslie2012}. 
In particular, observe that $r_k \ra 0$ since $\beta_k \ra 0$ and $\norm{B_\rho(p_k)}_\infty \leq \rho$ for all $k \geq 0$.

\bibliographystyle{apalike} 
\begin{small}
\bibliography{references}
\end{small}

\end{document}

%% file: AMIN_style.tex
\makeatletter
\def\@settitle{\begin{center}%
		\baselineskip14\p@\relax
		\normalfont\LARGE\bfseries
		\@title
	\end{center}%
}

\def\section{\@startsection{section}{1}%
	\z@{.7\linespacing\@plus\linespacing}{.5\linespacing}%
	{\normalfont\large\bfseries}}

\def\subsection{\@startsection{subsection}{2}%
	\z@{.5\linespacing\@plus.7\linespacing}{.5\linespacing}%
	{\normalfont\bfseries}}

\def\@setauthors{%
  \begingroup
  \def\thanks{\protect\thanks@warning}%
  \trivlist
  \centering\footnotesize \@topsep30\p@\relax
  \advance\@topsep by -\baselineskip
  \item\relax
  \author@andify\authors
  \def\\{\protect\linebreak}%
  \authors%
  \ifx\@empty\contribs
  \else
    ,\penalty-3 \space \@setcontribs
    \@closetoccontribs
  \fi
  \endtrivlist
  \endgroup
}

\makeatother

\usepackage[table, xcdraw, usenames, dvipsnames]{xcolor}
\definecolor{darkblue}{rgb}{0.0, 0.0, 0.45}
\definecolor{darkgreen}{rgb}{0.0, 0.45, 0}
\usepackage[colorlinks	= true,
raiselinks	= true,
linkcolor	= darkblue, 
citecolor	= Mahogany,
urlcolor	= darkgreen,
pdfauthor	= {},
pdftitle	= {},
pdfkeywords	= {},
pdfsubject	= {},
plainpages	= false]{hyperref}

\pdfoutput=1
\date{\today}

\usepackage{dsfont,amsfonts,amssymb,amsmath,amsthm}
\usepackage{mathtools, mathrsfs}

\usepackage{graphicx}
\usepackage[font=small,labelfont=rm]{subcaption}
\usepackage[font=small,margin=10pt]{caption}
\usepackage{wrapfig}

\usepackage{multirow}
\usepackage{booktabs}

\usepackage{algorithm,algorithmic}

\usepackage{enumitem}
\usepackage{tikz}
\usepackage{courier} 

\usepackage{nicefrac}

\theoremstyle{plain}
\newtheorem{Thm}{Theorem}[section]

\newtheorem{Lem}[Thm]{Lemma}

\newcommand{\Min}{\min\limits}

\DeclareMathOperator*{\argmin}{\arg\!\min}

\newcommand{\R}{\mathbb{R}}

\newcommand{\ra}{\rightarrow}

\newcommand{\wt}{\widetilde}
\newcommand{\wh}{\widehat}
\newcommand{\Let}{\coloneqq}

\newcommand{\tr}{^{\top}}

\newcommand{\norm}[1]{\left\Vert #1 \right\Vert}

\def\ssum{\begingroup\textstyle \sum\endgroup}

\newcommand{\opt}{^\star}
\newcommand{\PP}{\mathds{P}}
\newcommand{\EE}{\mathds{E}}

\newcommand{\set}[1]{\mathcal{#1}}

\DeclareMathOperator{\ord}{\mathcal{O}}

%% file: Tab/equivalence_table_4.tex
\input{Tab/colors.tex}

\begin{table}[]
\centering
\resizebox{\textwidth}{!}{%
\renewcommand{\arraystretch}{1.1}
{\arrayrulecolor{gray}
\begin{tabular}{cccc}
\hline
    \multicolumn{2}{c}{ \cellcolor{Opt} \textbf{Optimization}} 
  & \multicolumn{2}{c}{ \cellcolor{Con} \textbf{Control}}         \\ 
\hline
  \addlinespace[5pt]
\hline
    \multicolumn{2}{c}{ \cellcolor{Opt} \renewcommand{\arraystretch}{1.7} \begin{tabular}{c} $ \Min_x \{ f(x)\Let \EE_{\xi} [\hf(x,\xi)] \}$ \\ Function $\hf:\R^{\ell}\times\bf{\Xi}\ra\R$, \\Random variable $\xi \sim \PP$ \end{tabular}} 
  & \multicolumn{2}{c}{\cellcolor{Con} \renewcommand{\arraystretch}{1.7} \begin{tabular}{c} $\Min_{\pi:\set{S}\ra\set{A}}   \EE_{(s_t)_{t=0}^\infty} \left[ \sum_{t=0}^{\infty} \gamma^t c(s_t,a_t)\,|\, s_0 = s,\ a_t = \pi(s_t)\right] ,\ \forall s \in \set{S} $ \\ State $s \in \set{S}$, Control $a \in \set{A}$,  Cost $c:\set{S}\times\set{A}\ra\R$   \\ Dynamics $s^+ \sim \PP(\cdot|s,a)$ \end{tabular}} \\ 
  \hline
  \addlinespace[5pt]
\hline 
    \multicolumn{1}{c}{ \cellcolor{DetOpt} \textbf{Deterministic} } 
  & \multicolumn{1}{c}{ \cellcolor{StoOpt} \textbf{Stochastic} } 
  & \multicolumn{1}{c}{ \cellcolor{MBCon} \textbf{Model-based} } 
  & \multicolumn{1}{c}{ \cellcolor{MFCon} \textbf{Model-free} }  \\
\hline
  \addlinespace[2pt]
  \multicolumn{4}{c}{ \hspace{-40pt} \textbf{Equivalent characterization}} \\
  \hline
  
    \multicolumn{1}{c}{ \cellcolor{DetOpt} \renewcommand{\arraystretch}{1.5} \begin{tabular}{c} $\gr f (x\opt) = 0 $ \end{tabular}} 
  & \multicolumn{1}{c}{ \cellcolor{StoOpt} \renewcommand{\arraystretch}{1.5} \begin{tabular}{c} $ \EE_{\xi} [\gr \hf (x\opt,\xi)] = 0$  \end{tabular}} 
  & \multicolumn{1}{c}{ \cellcolor{MBCon} \renewcommand{\arraystretch}{1.5} \begin{tabular}{c} $v\opt = \bo (v\opt)$ \end{tabular}} 
  & \multicolumn{1}{c}{ \cellcolor{MFCon} \renewcommand{\arraystretch}{1.5} \begin{tabular}{c} $q\opt = \EE_{\varsigma^+} \big[\sbo (q\opt,\varsigma^+)  \big]$ \end{tabular}}  \\
\hline
  \addlinespace[2pt]
  \multicolumn{4}{c}{ \hspace{-40pt} \textbf{Available oracle/information}} \\
\hline
    \multicolumn{1}{c}{ \cellcolor{DetOpt} \renewcommand{\arraystretch}{1.5}  \begin{tabular}{c} $\gr f (x)$ \\ (Prob. distribution $\PP$) \\ \end{tabular}} 
  & \multicolumn{1}{c}{ \cellcolor{StoOpt} \renewcommand{\arraystretch}{1.5}  \begin{tabular}{c} $\gr \hf (x,\xi)$ \\ (Samples $\xi$)  \\ \end{tabular}} 
  & \multicolumn{1}{c}{ \cellcolor{MBCon}  \renewcommand{\arraystretch}{1.5} \begin{tabular}{c} $\bo (v)$ \\ (Prob. kernel $\PP$, Cost $c$) \\ \end{tabular}} 
  & \multicolumn{1}{c}{ \cellcolor{MFCon}  \renewcommand{\arraystretch}{1.5} \begin{tabular}{c} $\sbo (q,\varsigma^+)$ \\ (Samples $\{s,a,c(s,a),s^+\}$) \\ \end{tabular}}  \\
\hline
  \addlinespace[2pt]
  \multicolumn{4}{c}{ \hspace{-40pt} \textbf{Transformation}} \\
  \hline
  
    \multicolumn{1}{c}{ \cellcolor{DetOpt} \hspace{2cm} \begin{tabular}{r} $x$ \\ $\gr f$ \\ $\id - \gr f$ \\ $\gr^2 f$ \\ $I - \gr^2 f$ \\ \end{tabular}} 
  & \multicolumn{1}{c}{  \begin{tabular}{c} $\xleftrightarrow{\hspace*{3.5cm}}$\\ $\xleftrightarrow{\hspace*{3.5cm}}$\\ $\xleftrightarrow{\hspace*{3.5cm}}$  \\ $\xleftrightarrow{\hspace*{3.5cm}}$\\ $\xleftrightarrow{\hspace*{3.5cm}}$  \\ \end{tabular}} 
  & \multicolumn{1}{c}{ \cellcolor{MBCon} \hspace{-3.5cm}  \begin{tabular}{l} $v$ \\ $\id - \bo$ \\ $\bo $  \\ $I - \gamma P$ \\ $\gamma P $ \\ \end{tabular}} 
  & \multicolumn{1}{c}{  }  \\                           
    \multicolumn{1}{c}{  } 
  & \multicolumn{1}{c}{ \cellcolor{StoOpt} \hspace{3cm} \begin{tabular}{r} $(x,\xi)$ \\ $\gr \hf $ \\ $\id - \gr \hf $ \\ $\gr^2 \hf $ \\ $I - \gr^2 \hf $ \\ \end{tabular}} 
  & \multicolumn{1}{c}{  \begin{tabular}{c} $\xleftrightarrow{\hspace*{4cm}}$\\ $\xleftrightarrow{\hspace*{4cm}}$\\ $\xleftrightarrow{\hspace*{4cm}}$  \\ $\xleftrightarrow{\hspace*{4cm}}$\\ $\xleftrightarrow{\hspace*{4cm}}$  \\ \end{tabular}} 
  & \multicolumn{1}{c}{ \cellcolor{MFCon} \hspace{-3.7cm}  \begin{tabular}{l} $(q,\varsigma^+)$ \\ $\id - \sbo $ \\ $\sbo $  \\ $I - \gamma \wh{P}$ \\ $\gamma \wh{P} $  \\ \end{tabular}} \\
\hline

\end{tabular}%
}}
\caption{Equivalence transformations. See Section~\ref{sec:problem} for details. 
}
\label{tab:equivalence}
\end{table}


%% file: Tab/colors.tex
\colorlet{Title}{gray!10} 
\colorlet{BigTitle}{gray!30} 
\colorlet{Opt}{NavyBlue!15}
\colorlet{Con}{red!15}
\colorlet{DetOpt}{NavyBlue!12}
\colorlet{MBCon}{red!12}
\colorlet{StoOpt}{NavyBlue!18}
\colorlet{MFCon}{red!18} 

%% file: Tab/extended_algorithms_table.tex
\input{Tab/colors.tex}
\begin{table}[]
\resizebox{\textwidth}{!}{%
\centering
\renewcommand{\arraystretch}{1.1}
  {\arrayrulecolor{gray}
  \begin{tabular}{cc|cc}
    
    \hline
        \multicolumn{1}{c|}{ \cellcolor{DetOpt} \begin{tabular}{c} \textbf{Deterministic} \\ \textbf{optimization} \\   \end{tabular}} 
      & \multicolumn{1}{c|}{ \cellcolor{MBCon} \begin{tabular}{c} \textbf{Model-based} \\ \textbf{control} \\ \end{tabular}} 
      & \multicolumn{1}{c|}{ \cellcolor{StoOpt} \begin{tabular}{c} \textbf{Stochastic} \\ \textbf{optimization} \\ \end{tabular}} 
      & \multicolumn{1}{c}{ \cellcolor{MFCon} \begin{tabular}{c} \textbf{Model-free} \\ \textbf{control} \\ \end{tabular}} \\
    \hline
      \multicolumn{1}{c|}{ \begin{tabular}{c} ($y = x$) \\ $g(x) \Let \gr f (x)$ \\ $H(x) \Let \nabla^2 f (x)$ \\ \end{tabular}} 
    & \multicolumn{1}{c|}{ \begin{tabular}{c} ($y = v$) \\ $g(v) \Let v - \bo (v) $ \\ $H(v) \Let I - \gamma P(v) $  \\ \end{tabular}} 
    & \multicolumn{1}{c|}{ \begin{tabular}{c} ($y = x$) \\ $\wh{g}_k (x) \Let \gr \hf (x,\xi_k)$ \\ $\wh{H}_k (x) \Let \nabla^2 \hf (x,\xi_k)$ \\ \end{tabular}} 
    & \multicolumn{1}{c}{ \begin{tabular}{c} ($y = q$) \\ $\wh{g}_k (q) \Let q - \sbo (q,\varsigma^+_k) $ \\ $\wh{H}_k(q) \Let I - \gamma \wh{P} (q,\varsigma^+_k) $ \\ \end{tabular}}  \\ 
    
    \hline
      \addlinespace[6pt]
    \hline
      \multicolumn{4}{c}{ \hspace{-45pt} \cellcolor{BigTitle} \textbf{First Order Methods (access to gradient oracle)}} \\
    \hline
      \multicolumn{4}{c}{ \hspace{-45pt} \cellcolor{Title} \textbf{Vanilla First Order Method}} \\
    \hline
    
      \multicolumn{1}{c}{ \begin{tabular}{c} GD~\cite{lemarechal2012cauchy} \\ \end{tabular}}
    & \multicolumn{1}{c|}{ \begin{tabular}{c} Rel-VI~\cite{bellman1957markovian,kushner1971accelerated} \\ \end{tabular}}
    & \multicolumn{1}{c}{ \begin{tabular}{c} SGD~\cite{robbins1951stochastic} \\ \end{tabular}}
    & \multicolumn{1}{c}{ \begin{tabular}{c} QL~\cite{watkins1992q} \\ \end{tabular}} \\

      \multicolumn{2}{c|}{ \begin{tabular}{c} $d_k = - \alpha_k g(y_k)$ \\ \end{tabular}} 
    & \multicolumn{2}{c}{ \begin{tabular}{c} $d_k = - \alpha_k \wh{g}_k(y_k)$  \\ \end{tabular}} \\
    
    \hline
        \multicolumn{4}{c}{ \hspace{-45pt} \cellcolor{Title} \textbf{Acceleration via Momentum}} \\
    \hline

      \multicolumn{1}{c}{ \begin{tabular}{c} Mom-GD \cite{polyak1964some} \\ \end{tabular} } 
    & \multicolumn{1}{c|}{ \begin{tabular}{c} Mom-VI \cite{goyal2019first} \\ \end{tabular} } 
    & \multicolumn{1}{c}{ \multirow{2}{*}{ \begin{tabular}{c} Mom-SGD \cite{yang2016unified} \\ \end{tabular} } }
    & \multicolumn{1}{c}{ \multirow{2}{*}{ \begin{tabular}{c} SQL \cite{ghavamzadeh2011speedy}, NeSA \cite{devraj2019matrix}, \\ Mom-QL \cite{weng2020momentum} \\ \end{tabular} } } \\ 
      
      \multicolumn{2}{c|}{ \begin{tabular}{c} $d_k = - \alpha_k g(y_k) + \beta_k d_{k-1}$ \\ \end{tabular} }  
    & \multicolumn{2}{c}{ }  \\
      \cline{1-2}
      
      \multicolumn{1}{c}{ \begin{tabular}{c} Acc-GD \cite{nesterov1983method} \\ \end{tabular} } 
    & \multicolumn{1}{c|}{ \begin{tabular}{c} Acc-VI \cite{goyal2019first} \\ \end{tabular} } 
    & \multicolumn{2}{c}{ \multirow{2}{*}{ \begin{tabular}{c} $\left\{\begin{array}{l} d'_{k-1} =  \wh{g}_k(y_{k}) - \wh{g}_k(y_{k-1}) \\ d_k =  -  \alpha_k \wh{g}_k(y_k) - \beta_k d'_{k-1}  + \delta_k d_{k-1} \end{array}\right.$\\ \end{tabular} } } \\
      \multicolumn{2}{c|}{ \begin{tabular}{c} $d_k = - \alpha_k g(y_k+\beta_kd_{k-1}) + \beta_k d_{k-1}$ \\ \end{tabular} } 
    &  \multicolumn{2}{c}{ } \\
  
    \hline
        \multicolumn{4}{c}{ \hspace{-45pt} \cellcolor{Title} \textbf{Acceleration via Anchoring}} \\
    \hline
  
      \multicolumn{1}{c}{ \begin{tabular}{c} Anc-GD~\cite{yoon2021accelerated} \\ \end{tabular}} 
    & \multicolumn{1}{c|}{ \begin{tabular}{c} Anc-VI~\cite{lee2023accelerating} \\ \end{tabular}} 
    & \multicolumn{1}{c}{ \begin{tabular}{c} SHI~\cite{cai2022stochastic} \\ \end{tabular}} 
    & \multicolumn{1}{c}{ \begin{tabular}{c} HQL~\cite{bravo2024stochastic}, SAVID~\cite{lee2025near} \\ \end{tabular}}  \\ 
  
      \multicolumn{2}{c|}{ \begin{tabular}{c} $d_k = \beta_k(y_0 - y_k) -\alpha_k g(y_k) $ \end{tabular}} 
    & \multicolumn{2}{c}{ \begin{tabular}{c} $d_k = \beta_k(y_0 - y_k) + \alpha_k \wh{g}_k(y_k)$  \end{tabular}} \\
  
    \hline
      \multicolumn{4}{c}{ \hspace{-45pt} \cellcolor{Title} \textbf{Acceleration via PID Control}} \\
    \hline
  
      \multicolumn{1}{c}{ \begin{tabular}{c} PID-Opt (Full-batch) \cite{an2018pid} \\ \end{tabular}} 
    & \multicolumn{1}{c|}{ \begin{tabular}{c} PID-VI \cite{farahmand2021pid} \\ \end{tabular}} 
    & \multicolumn{1}{c}{ \begin{tabular}{c} PID-Opt (Mini-batch)  \cite{an2018pid} \\ \end{tabular}} 
    & \multicolumn{1}{c}{ \begin{tabular}{c} PID-QL \cite{bedaywi2024pid} \\ \end{tabular}}  \\ 
  
      \multicolumn{2}{c|}{ \begin{tabular}{c} $\left\{\begin{array}{l}
        d^I_{k} = - \alpha_k g(y_k) + \beta_k d^I_{k-1}  \\
        d'_{k-1} = g(y_k) - g(y_{k-1}) \\
        d^D_{k} = \delta_k d_{k-1} + (1-\delta_k) d'_{k-1} \\
        d_k = - \kappa^P_k g(y_k) + \kappa^I_k d^I_{k}  + \kappa^D_k d^D_{k}
        \end{array}\right.$ \end{tabular}} 
    & \multicolumn{2}{c}{ \begin{tabular}{c} $\left\{\begin{array}{l}
        d^I_{k} = - \alpha_k \wh{g}_{k} (y_{k}) + \beta_k d^I_{k-1} \\
        d'_{k-1} = \wh{g}_{k}(y_k) - \wh{g}_{k-1}(y_{k-1}) \\
        y'_{k-  1} = (1-\eta_k) y'_{k-2} + \eta_k y_{k-1} \\
        d^D_{k} =  \delta_k (y_{k} - y'_{k-1}) + (1-\delta_k) d'_{k-1} \\
        d_k = - \kappa^P_k \wh{g}_{k} (y_{k}) + \kappa^I_k d^I_{k+1}  + \kappa^d_k d^D_{k}
      \end{array}\right.$  \end{tabular}} \\

    \hline
      \multicolumn{4}{c}{ \hspace{-45pt} \cellcolor{Title} \textbf{QNM via Anderson Mixing}$^\dagger$} \\
    \hline
  
      \multicolumn{1}{c}{ \begin{tabular}{c} AA-GD \cite{mai2020anderson} \\ \end{tabular}} 
    & \multicolumn{1}{c|}{ \begin{tabular}{c} AA-VI \cite{geist2018anderson} \\ \end{tabular}} 
    & \multicolumn{1}{c}{ \begin{tabular}{c} SAM \cite{NEURIPS2021_c203e4a1} \\ \end{tabular}} 
    & \multicolumn{1}{c}{ \begin{tabular}{c} SAA-QL~\cite{sun2021damped} \end{tabular}}  \\ 
  
      \multicolumn{2}{c|}{ \begin{tabular}{c} $\left\{\begin{array}{l}
             Y_k = [y_{k}, \ldots, y_{k-m_k}] \\
             G_k = [g(y_{k}), \ldots, g(y_{k-m_k})]\\
             w_k = (\e^\top (G_k^\top G_k)^{-1} \e)^{-1} (G_k^\top G_k)^{-1} \e\\
             d_k = - y_k + (Y_k - \alpha_k G_k)w_k
             \end{array}\right.$ \end{tabular}} 
    & \multicolumn{2}{c}{ \begin{tabular}{c} $\left\{\begin{array}{l}
             D^{y}_k = [d^{y}_{k-m_k}, \ldots, d^{y}_{k-1}] \quad \text{with} \quad d^{y}_{i} = y_{i+1} - y_{i} \\
             D^{g}_k = [d^{g}_{k-m_k}, \ldots, d^{g}_{k-1}] \quad \text{with} \quad d^{g}_{i} = \wh{g}_{i+1} (y_{i+1}) - \wh{g}_{i} (y_{i})\\
            D_k = \beta_k I - \alpha_k(D^{y}_k + \beta_k D^{g}_k)\big((D^{g}_k)^\top D^{g}_k + R_k\big)^{-1}(D^{g}_k)^\top\\
            d_k = - D_k \wh{g}_{k} (y_{k})
            \end{array}\right.$ \end{tabular}} \\
  
  \hline
    \multicolumn{4}{c}{ \hspace{-45pt} \cellcolor{BigTitle} \textbf{Second-Order Methods (access to gradient and Hessian oracles)}} \\
  
  \hline
    \multicolumn{4}{c}{ \hspace{-45pt} \cellcolor{Title} \textbf{QNM via Spectral Preconditioning}$^\ddagger$} \\
  \hline
  
    \multicolumn{1}{c}{ \begin{tabular}{c} (SP-GD)\cite{doikov2024spectral} \\ \end{tabular}} 
  & \multicolumn{1}{c|}{ \begin{tabular}{c} R1-VI \cite{kolarijani2025rank} \end{tabular}} 
  & \multicolumn{1}{c}{ \begin{tabular}{c} (NewSamp)~\cite{erdogdu2015convergence} \\ \end{tabular}} 
  & \multicolumn{1}{c}{ \begin{tabular}{c} R1-QL \cite{kolarijani2025rank} \\ \end{tabular}}  \\ 
  
      \multicolumn{2}{c|}{ \begin{tabular}{c} $\left\{\begin{array}{l}
         D_k = \alpha_k I + \beta_k U_k W_k\tr ,\\
         d_k = - D_k^{-1} g(y_k),
         \end{array}\right.$ \end{tabular}} 
    & \multicolumn{2}{c}{ \begin{tabular}{c} $\left\{\begin{array}{l}
          D_k = \alpha_k I + \beta_k \wh{U}_k \wh{W}_k\tr ,\\
          d_k = - D_k^{-1} \wh{g}_k(y_k)
          \end{array}\right.$  \end{tabular}} \\
  
  \hline
    \multicolumn{4}{c}{ \hspace{-45pt} \cellcolor{Title} \textbf{Vanilla Second-Order Method}} \\
  \hline
  
      \multicolumn{1}{c}{ \begin{tabular}{c} NM \\ \end{tabular}} 
    & \multicolumn{1}{c|}{ \begin{tabular}{c} PI \cite{howard1960dynamic} \\ \end{tabular}} 
    & \multicolumn{1}{c}{ \begin{tabular}{c} SNR \cite{ruppert1985newton,DevrajMeyn2017NIPS} \\ \end{tabular}} 
    & \multicolumn{1}{c}{ \begin{tabular}{c} ZQL \cite{DevrajMeyn2017NIPS} \\ \end{tabular}}  \\ 
  
      \multicolumn{2}{c|}{ \begin{tabular}{c} $d_k = - \alpha_k \left[ H(y_k) \right]^{-1} g(y_{k}) $ \\ \end{tabular}} 
    & \multicolumn{2}{c}{ \begin{tabular}{c} $\left\{\begin{array}{l} D_k = (1-\beta_k) D_{k-1} + \beta_k \wh{H}_k(y_k) \\ d_k = -\alpha_k D_k^{-1} \wh{g}_k(y_k) \end{array}\right.$  \\ \end{tabular}} \\
    
    \toprule
  \end{tabular}%
}}
\caption{Equivalent algorithms: The vector~$d_k$ is the update vector in the iterative scheme $y_{k+1} = y_{k}+d_k$ for $k= 0,1\ldots$. The Greek letters $\alpha_k, \beta_k, \ldots$ are scalar coefficients that may take different values (even in equivalent algorithms). $^\dagger$$R_k$ is a regularizer depending on $D^y_k$ and $D^g_k$. $^\ddagger$$U_k$ and $W_k$ are tall matrices depending on the (sampled) second-order oracle. See Section~\ref{sec:equiv_alg} for details.
}
\label{tab:algorithms}
\end{table}

%% file: references.bib
@book{bertsekas2022lessons,
  title={Lessons from AlphaZero for optimal, model predictive, and adaptive control},
  author={Bertsekas, Dimitri},
  year={2022},
  publisher={Athena Scientific}
}

@book{villani2008optimal,
  title={Optimal transport: old and new},
  author={Villani, C{\'e}dric and others},
  volume={338},
  year={2008},
  publisher={Springer}
}

@inproceedings{erdogdu2015convergence,
	title = {Convergence rates of sub-sampled {N}ewton methods},
	author = {Erdogdu, Murat A and Montanari, Andrea},
	booktitle = {Advances in Neural Information Processing Systems},
	volume = {28},
	year = {2015},
}

@inproceedings{doikov2024spectral,
	title = {Spectral Preconditioning for Gradient Methods on Graded
	         Non-convex Functions},
	author = {Doikov, Nikita and Stich, Sebastian U and Jaggi, Martin},
	booktitle = {International Conference on Machine Learning},
	pages = {11227--11252},
	year = {2024},
}

@article{liu2024sketched,
	title = {Sketched {N}ewton value iteration for large-scale {M}arkov
	         decision processes},
	author = {Liu, Jinsong and Xie, Chenghan and Deng, Qi and Ge, Dongdong
	          and Ye, Yinyu},
	journal = {Proceedings of the AAAI Conference on Artificial Intelligence},
	volume = {38},
	number = {12},
	pages = {13936--13944},
	year = {2024},
}

@article{fang2009two,
	title = {Two classes of multisecant methods for nonlinear acceleration},
	author = {Fang, Haw-ren and Saad, Yousef},
	journal = {Numerical Linear Algebra with Applications},
	volume = {16},
	number = {3},
	pages = {197--221},
	year = {2009},
	publisher = {Wiley Online Library},
}

@article{anderson1965iterative,
	title = {Iterative procedures for nonlinear integral equations},
	author = {Anderson, Donald G},
	journal = {Journal of the ACM (JACM)},
	volume = {12},
	number = {4},
	pages = {547--560},
	year = {1965},
	publisher = {ACM New York, NY, USA},
}

@article{kolarijani2023optimization,
	title = {From optimization to control: quasi policy iteration},
	author = {Kolarijani, Mohammad Amin Sharifi and Mohajerin Esfahani,
	          Peyman},
	journal = {arXiv preprint arXiv:2311.11166},
	year = {2025},
}

@inproceedings{jin2024truncated,
	title = {Truncated variance reduced value iteration},
	author = {Jin, Yujia and Karmarkar, Ishani and Sidford, Aaron and Wang,
	          Jiayi},
	booktitle = {Advances in Neural Information Processing Systems},
	volume = {37},
	pages = {117481--117508},
	year = {2024},
}

@inproceedings{park2022exact,
	title = {Exact optimal accelerated complexity for fixed-point iterations
	         },
	author = {Park, Jisun and Ryu, Ernest K.},
	booktitle = {International Conference on Machine Learning},
	pages = {17420--17457},
	year = {2022},
}

@article{fletcher1963rapidly,
	title = {A rapidly convergent descent method for minimization},
	author = {Fletcher, Roger and Powell, Michael JD},
	journal = {Computer Journal},
	volume = {6},
	number = {2},
	pages = {163--168},
	year = {1963},
	publisher = {The British Computer Society},
}

@article{davidon1991variable,
	title = {Variable metric method for minimization},
	author = {Davidon, William C},
	journal = {SIAM Journal on Optimization},
	volume = {1},
	number = {1},
	pages = {1--17},
	year = {1991},
	publisher = {SIAM},
}

@incollection{powell1970new,
	title = {A new algorithm for unconstrained optimization},
	author = {Powell, Michael JD},
	booktitle = {Nonlinear programming},
	pages = {31--65},
	year = {1970},
	publisher = {Elsevier},
}

@article{Liu2024fitted,
	author = {Liu, Yichen and Kolarijani, Mohamad Amin Sharifi },
	journal = {IEEE Control Systems Letters},
	title = {Fitted {Q}-iteration via {M}ax-{P}lus-linear approximation},
	year = {2024},
	volume = {8},
	pages = {3201--3206},
	publisher = {IEEE},
}

@article{nesterov2005smooth,
	title = {Smooth minimization of non-smooth functions},
	author = {Nesterov, Yurii},
	journal = {Mathematical Programming},
	volume = {103},
	number = {1},
	pages = {127--152},
	year = {2005},
	publisher = {Springer},
}

@article{halpern1967fixed,
	title = {Fixed points of nonexpanding maps},
	author = {Halpern, Benjamin},
	journal = {Bulletin of the American Mathematical Society},
	volume = {73},
	number = {6},
	pages = {957--961},
	year = {1967},
	publisher = {AMS},
}

@inproceedings{sun2021damped,
	title = {Damped {A}nderson mixing for deep reinforcement learning:
	         acceleration, convergence, and stabilization},
	author = {Sun, Ke and Wang, Yafei and Liu, Yi and Pan, Bo and Jui,
	          Shangling and Jiang, Bei and Kong, Linglong and others},
	booktitle = {Advances in Neural Information Processing Systems},
	volume = {34},
	pages = {3732--3743},
	year = {2021},
}

@book{szepesvari2022algorithms,
	title = {Algorithms for reinforcement learning},
	author = {Szepesv{\'a}ri, Csaba},
	year = {2010},
	publisher = {Morgan \& Claypool},
}

@article{GONCALVES2021109623,
	title = {{M}ax-{P}lus approximation for reinforcement learning},
	journal = {Automatica},
	volume = {129},
	pages = {109623},
	year = {2021},
	author = {Vinicius Mariano Gon\c{c}alves},
	publisher = {Elsevier},
}

@book{McEn06,
	title = {{M}ax-{P}lus methods for nonlinear control and estimation},
	author = {McEneaney, William M},
	year = {2006},
	publisher = {Springer Science \& Business Media},
}

@article{Bach20,
	author = {Berthier, Eloïse and Bach, Francis},
	journal = {IEEE Control Systems Letters},
	title = {{M}ax-{P}lus linear approximations for deterministic
	         continuous-state {M}arkov decision processes},
	year = {2020},
	volume = {4},
	number = {3},
	pages = {767--772},
	publisher = {IEEE},
}

@article{Tsitsiklis97,
	author = {Tsitsiklis, John N. and Van Roy, Benjamin},
	journal = {IEEE Transactions on Automatic Control},
	title = {An analysis of temporal-difference learning with function
	         approximation},
	year = {1997},
	volume = {42},
	number = {5},
	pages = {674--690},
	publisher = {IEEE},
}

@article{SCHMIDHUBER2015,
	title = {Deep learning in neural networks: an overview},
	journal = {Neural Networks},
	volume = {61},
	pages = {85--117},
	year = {2015},
	author = {J{\"u}rgen Schmidhuber},
	publisher = {Elsevier},
}

@book{sutton2018reinforcement,
	title = {Reinforcement learning: an introduction},
	author = {Sutton, Richard S. and Barto, Andrew G.},
	year = {2018},
	publisher = {MIT press},
}

@article{broyden1965class,
	title = {A class of methods for solving nonlinear simultaneous equations
	         },
	author = {Broyden, Charles G},
	journal = {Mathematics of Computation},
	volume = {19},
	number = {92},
	pages = {577--593},
	year = {1965},
	publisher = {AMS},
}

@article{bellman1957markovian,
	title = {A {M}arkovian decision process},
	author = {Bellman, Richard},
	journal = {Journal of Mathematics and Mechanics},
	volume = {6},
	number = {5},
	pages = {679--684},
	year = {1957},
	publisher = {JSTOR},
}

@article{goyal2019first,
	title = {A first-order approach to accelerated value iteration},
	author = {Goyal, Vineet and Grand-Cl{\'e}ment, Julien},
	journal = {Operations Research},
	volume = {71},
	number = {2},
	pages = {517--535},
	year = {2022},
	publisher = {INFORMS},
}

@inproceedings{ghavamzadeh2011speedy,
	title = {Speedy {Q}-learning},
	author = {Ghavamzadeh, Mohammad and Kappen, Hilbert and Azar, Mohammad
	          and Munos, R{\'e}mi},
	booktitle = {Advances in Neural Information Processing Systems},
	volume = {24},
	year = {2011},
}

@inproceedings{weng2020momentum,
	title = {Finite-time theory for momentum {Q}-learning},
	author = {Weng, Bowen and Xiong, Huaqing and Zhao, Lin and Liang,
	          Yingbin and Zhang, Wei},
	booktitle = {Conference on Uncertainty in Artificial Intelligence},
	pages = {665--674},
	year = {2021},
}

@article{polyak1964some,
	title = {Some methods of speeding up the convergence of iteration
	         methods},
	author = {Polyak, Boris T},
	journal = {{USSR} Computational Mathematics and Mathematical Physics},
	volume = {4},
	number = {5},
	pages = {1--17},
	year = {1964},
	publisher = {Elsevier},
}

@inproceedings{nesterov1983method,
	title = {A method for solving the convex programming problem with
	         convergence rate $O (1/k^2)$},
	author = {Nesterov, Yurii},
	booktitle = {Doklady Akademii Nauk SSSR},
	volume = {269},
	pages = {543--547},
	year = {1983},
}

@article{watkins1992q,
	title = {{Q}-learning},
	author = {Watkins, Christopher JCH and Dayan, Peter},
	journal = {Machine Learning},
	volume = {8},
	number = {3},
	pages = {279--292},
	year = {1992},
	publisher = {Springer},
}

@article{lemarechal2012cauchy,
	title = {{C}auchy and the gradient method},
	author = {Lemar{\'e}chal, Claude},
	journal = {Documenta Mathematica Extra},
	pages = {251--254},
	year = {2012},
	publisher = {Deutsche Mathematiker-Vereinigung},
}

@article{kamanchi2021generalized,
	title = {Generalized second order value iteration in {M}arkov decision
	         processes},
	author = {Kamanchi, Chandramouli and Diddigi, Raghuram Bharadwaj and
	          Bhatnagar, Shalabh},
	journal = {IEEE Transactions on Automatic Control},
	volume = {67},
	number = {8},
	pages = {4241--4247},
	year = {2022},
	publisher = {IEEE},
}

@article{puterman1979convergence,
	title = {On the convergence of policy iteration in stationary dynamic
	         programming},
	author = {Puterman, Martin L and Brumelle, Shelby L},
	journal = {Mathematics of Operations Research},
	volume = {4},
	number = {1},
	pages = {60--69},
	year = {1979},
	publisher = {INFORMS},
}

@article{porteus1978accelerated,
	title = {Accelerated computation of the expected discounted return in a
	         {M}arkov chain},
	author = {Porteus, Evan L and Totten, John C},
	journal = {Operations Research},
	volume = {26},
	number = {2},
	pages = {350--358},
	year = {1978},
	publisher = {INFORMS},
}

@article{robbins1951stochastic,
	title = {A stochastic approximation method},
	author = {Robbins, Herbert and Monro, Sutton},
	journal = {Annals of Mathematical Statistics},
	volume = {22},
	number = {3},
	pages = {400--407},
	year = {1951},
	publisher = {JSTOR},
}

@article{ruppert1985newton,
	title = {A {N}ewton-{R}aphson version of the multivariate {R}obbins-{M}
	         onro procedure},
	author = {Ruppert, David},
	journal = {Annals of Statistics},
	volume = {13},
	number = {1},
	pages = {236--245},
	year = {1985},
	publisher = {Institute of Mathematical Statistics},
}

@inproceedings{devraj2019matrix,
	title = {On matrix momentum stochastic approximation and applications to
	         {Q}-learning},
	author = {Devraj, Adithya M. and Bu{\v{s}}i{\'c}, Ana and Meyn, Sean P.},
	booktitle = {Annual Allerton Conference on Communication, Control, and
	             Computing},
	pages = {749--756},
	year = {2019},
}

@book{howard1960dynamic,
	title = {Dynamic programming and {M}arkov processes},
	author = {Howard, Ronald A},
	year = {1960},
	publisher = {John Wiley},
}

@article{kushner1971accelerated,
	title = {Accelerated procedures for the solution of discrete {M}arkov
	         control problems},
	author = {Kushner, Harold and Kleinman, A},
	journal = {IEEE Transactions on Automatic Control},
	volume = {16},
	number = {2},
	pages = {147--152},
	year = {1971},
	publisher = {IEEE},
}

@inproceedings{liu2018accelerating,
	title = {Accelerating {SGD} with momentum for over-parameterized
	         learning},
	author = {Liu, Chaoyue and Belkin, Mikhail},
	booktitle = {International Conference on Learning Representations},
	year = {2020},
}

@inproceedings{kidambi2018insufficiency,
	title = {On the insufficiency of existing momentum schemes for
	         stochastic optimization},
	author = {Kidambi, Rahul and Netrapalli, Praneeth and Jain, Prateek and
	          Kakade, Sham},
	booktitle = {Information Theory and Applications Workshop},
	pages = {1--9},
	year = {2018},
}

@article{samuelson1948foundations,
	title = {Foundations of economic analysis},
	author = {Samuelson, Paul Anthony},
	journal = {Science and Society},
	volume = {13},
	number = {1},
	year = {1948},
	publisher = {Guilford Press},
}

@article{grand2021convex,
	title = {From convex optimization to {MDP}s: a review of first-order,
	         second-order and quasi-{N}ewton methods for {MDP}s},
	author = {Grand-Cl{\'e}ment, Julien},
	journal = {arXiv preprint arXiv:2104.10677},
	year = {2021},
}

@article{bubeck2015convex,
	title = {Convex optimization: algorithms and complexity},
	author = {Bubeck, S{\'e}bastien},
	journal = {Foundations and Trends{\textregistered} in Machine Learning},
	volume = {8},
	number = {3-4},
	pages = {231--357},
	year = {2015},
	publisher = {Now Publishers, Inc.},
}

@inproceedings{liu2022almost,
	title = {On almost sure convergence rates of stochastic gradient methods
	         },
	author = {Liu, Jun and Yuan, Ye},
	booktitle = {Conference on Learning Theory},
	pages = {2963--2983},
	year = {2022},
}

@book{nesterov2018lectures,
	title = {Lectures on convex optimization},
	author = {Nesterov, Yurii},
	year = {2018},
	publisher = {Springer},
}

@article{allen2017katyusha,
	title = {Katyusha: the first direct acceleration of stochastic gradient
	         methods},
	author = {Allen-Zhu, Zeyuan},
	journal = {Journal of Machine Learning Research},
	volume = {18},
	number = {1},
	pages = {8194--8244},
	year = {2017},
	publisher = {JMLR. org},
}

@inproceedings{kearns1998finite,
	title = {Finite-sample convergence rates for {Q}-learning and indirect
	         algorithms},
	author = {Kearns, Michael and Singh, Satinder},
	booktitle = {Advances in Neural Information Processing Systems},
	volume = {11},
	year = {1998},
}

@article{rust1994structural,
	title = {Structural estimation of {M}arkov decision processes},
	author = {Rust, John},
	journal = {Handbook of Econometrics},
	volume = {4},
	pages = {3081--3143},
	year = {1994},
	publisher = {Elsevier},
}

@article{yang2016unified,
	title = {Unified convergence analysis of stochastic momentum methods for
	         convex and non-convex optimization},
	author = {Yang, T. and Lin, Q. and Li, Z.},
	journal = {arXiv preprint arXiv:1604.03257},
	year = {2016},
}

@article{zhang2020globally,
	title = {Globally convergent type-{I} {A}nderson acceleration for
	         nonsmooth fixed-point iterations},
	author = {Zhang, Junzi and O'Donoghue, Brendan and Boyd, Stephen},
	journal = {SIAM Journal on Optimization},
	volume = {30},
	number = {4},
	pages = {3170--3197},
	year = {2020},
	publisher = {SIAM},
}

@article{geist2018anderson,
	title = {{A}nderson acceleration for reinforcement learning},
	author = {Geist, Matthieu and Scherrer, Bruno},
	journal = {arXiv preprint arXiv:1809.09501},
	year = {2018},
}

@article{gargiani2022dynamic,
	title = {Dynamic programming through the lens of semismooth {N}
	         ewton-type methods},
	author = {Gargiani, Matilde and Zanelli, Andrea and Liao-McPherson,
	          Dominic and Summers, TH and Lygeros, John},
	journal = {IEEE Control Systems Letters},
	volume = {6},
	pages = {2996--3001},
	year = {2022},
	publisher = {IEEE},
}

@article{vieillard2019connections,
	title = {On connections between constrained optimization and
	         reinforcement learning},
	author = {Vieillard, Nino and Pietquin, Olivier and Geist, Matthieu},
	journal = {arXiv preprint arXiv:1910.08476},
	year = {2019},
}

@inproceedings{kakade2002approximately,
	title = {Approximately optimal approximate reinforcement learning},
	author = {Kakade, Sham and Langford, John},
	booktitle = {International Conference on Machine Learning},
	pages = {267--274},
	year = {2002},
}

@article{frank1956algorithm,
	title = {An algorithm for quadratic programming},
	author = {Frank, Marguerite and Wolfe, Philip},
	journal = {Naval Research Logistics Quarterly},
	volume = {3},
	number = {1-2},
	pages = {95--110},
	year = {1956},
	publisher = {Wiley},
}

@book{ref:Bert_neuro-dynamic,
	title = {Neuro-dynamic programming},
	author = {Bertsekas, Dimitri P. and Tsitsiklis, John N.},
	year = {1996},
	publisher = {Athena Scientific},
}

@book{ref:HL1,
	title = {Discrete-time {M}arkov control processes: basic optimality
	         criteria},
	author = {Hern{\'a}ndez-Lerma, On{\'e}simo and Lasserre, Jean B},
	volume = {30},
	year = {2012},
	publisher = {Springer Science \& Business Media},
}

@book{ref:HL2,
	title = {Further topics on discrete-time {M}arkov control processes},
	author = {Hern{\'a}ndez-Lerma, On{\'e}simo and Lasserre, Jean B},
	volume = {42},
	year = {2012},
	publisher = {Springer Science \& Business Media},
}

@article{ref:Bert_disc_75,
	title = {Convergence of discretization procedures in dynamic programming
	         },
	author = {Bertsekas, Dimitri P.},
	journal = {IEEE Transactions on Automatic Control},
	volume = {20},
	number = {3},
	pages = {415--419},
	year = {1975},
	publisher = {IEEE},
}

@book{ref:Powell_07,
	title = {Approximate dynamic programming: solving the curses of
	         dimensionality},
	author = {Powell, Warren B},
	volume = {703},
	year = {2007},
	publisher = {John Wiley \& Sons},
}

@article{ref:Bert_temporal,
	title = {Temporal difference methods for general projected equations},
	author = {Bertsekas, Dimitri P.},
	journal = {IEEE Transactions on Automatic Control},
	volume = {56},
	number = {9},
	pages = {2128--2139},
	year = {2011},
	publisher = {IEEE},
}

@book{ref:Bert_abstract,
	title = {Abstract dynamic programming},
	author = {Bertsekas, Dimitri P.},
	year = {2022},
	publisher = {Athena Scientific},
}

@article{ref:FDP_TAC,
	author = {Kolarijani, Mohamad Amin Sharifi and {Mohajerin Esfahani},
	          Peyman },
	journal = {IEEE Transactions on Automatic Control},
	title = {Fast approximate dynamic programming for input-affine dynamics},
	year = {2023},
	volume = {68},
	number = {10},
	pages = {6315--6322},
	publisher = {IEEE},
}

@inproceedings{ref:FDP_NeurIPS,
	title = {Fast approximate dynamic programming for infinite-horizon {M}
	         arkov decision processes},
	author = {Kolarijani, Mohamad Amin Sharifi and Max, Gyula F and {
	          Mohajerin Esfahani}, Peyman},
	booktitle = {Advances in Neural Information Processing Systems},
	volume = {34},
	pages = {23652--23663},
	year = {2021},
}

@article{ref:Moh_SIOPT,
	title = {From infinite to finite programs: explicit error bounds with
	         applications to approximate dynamic programming},
	author = {{Mohajerin Esfahani}, Peyman and Sutter, Tobias and Kuhn,
	          Daniel and Lygeros, John},
	journal = {SIAM Journal on Optimization},
	volume = {28},
	number = {3},
	pages = {1968--1998},
	year = {2018},
	publisher = {SIAM},
}

@article{ref:VanRoy_MP,
	title = {On constraint sampling in the linear programming approach to
	         approximate dynamic programming},
	author = {De Farias, Daniela Pucci and Van Roy, Benjamin},
	journal = {Mathematics of Operations Research},
	volume = {29},
	number = {3},
	pages = {462--478},
	year = {2004},
	publisher = {INFORMS},
}

@inproceedings{shi2019regularized,
	title = {Regularized {A}nderson acceleration for off-policy deep
	         reinforcement learning},
	author = {Shi, Wenjie and Song, Shiji and Wu, Hui and Hsu, Ya-Chu and Wu
	          , Cheng and Huang, Gao},
	booktitle = {Advances in Neural Information Processing Systems},
	volume = {32},
	year = {2019},
}

@inproceedings{lee2023accelerating,
	title = {Accelerating value iteration with anchoring},
	author = {Lee, Jongmin and Ryu, Ernest K.},
	booktitle = {Advances in Neural Information Processing Systems},
	volume = {36},
	pages = {53924--53963},
	year = {2023},
}

@inproceedings{farahmand2021pid,
	title = {{PID} accelerated value iteration algorithm},
	author = {Farahmand, Amir-Massoud and Ghavamzadeh, Mohammad},
	booktitle = {International Conference on Machine Learning},
	pages = {3143--3153},
	year = {2021},
}

@article{bedaywi2024pid,
	title = {{PID} accelerated temporal difference algorithms},
	author = {Bedaywi, Mark and Rakhsha, Amin and Farahmand, Amir-massoud},
	journal = {Reinforcement Learning Journal},
	volume = {5},
	pages = {2071--2095},
	year = {2025},
	publisher = {OpenReview},
}

@article{lee2024deflated,
	title = {Deflated dynamics value iteration},
	author = {Lee, Jongmin and Rakhsha, Amin and Ryu, Ernest K. and
	          Farahmand, Amir-massoud},
	journal = {Transactions on Machine Learning Research},
	year = {2025},
	publisher = {OpenReview},
}

@inproceedings{kolarijani2025rank,
	title = {Rank-one modified value iteration},
	author = {Kolarijani, Arman Sharifi and Ok, Tolga and Mohajerin Esfahani
	          , Peyman and Sharifi Kolarijani, Mohamad Amin},
	booktitle = {International Conference on Machine Learning},
	pages = {31182--31201},
	year = {2025},
}

@inproceedings{NEURIPS2021_c203e4a1,
	author = {Wei, Fuchao and Bao, Chenglong and Liu, Yang},
	booktitle = {Advances in Neural Information Processing Systems},
	pages = {22995--23008},
	title = {Stochastic {A}nderson mixing for nonconvex stochastic
	         optimization},
	volume = {34},
	year = {2021},
}

@inproceedings{mai2020anderson,
	title = {{A}nderson acceleration of proximal gradient methods},
	author = {Mai, Vien and Johansson, Mikael},
	booktitle = {International Conference on Machine Learning},
	pages = {6620--6629},
	year = {2020},
}

@article{zuo2022offline,
	title = {Offline reinforcement learning with {A}nderson acceleration for
	         robotic tasks},
	author = {Zuo, Guoyu and Huang, Shuai and Li, Jiangeng and Gong,
	          Daoxiong},
	journal = {Applied Intelligence},
	volume = {52},
	number = {9},
	pages = {9885--9898},
	year = {2022},
	publisher = {Springer},
}

@article{tran2022connection,
	title = {From {H}alpern's fixed-point iterations to {N}esterov's
	         accelerated interpretations for root-finding problems},
	author = {Tran-Dinh, Quoc},
	journal = {Computational Optimization and Applications},
	volume = {87},
	number = {1},
	pages = {181--218},
	year = {2024},
	publisher = {Springer},
}

@inproceedings{yoon2021accelerated,
	title = {Accelerated algorithms for smooth convex-concave minimax
	         problems with $\mathcal{O}(1/k^2)$ rate on squared gradient norm
	         },
	author = {Yoon, Taeho and Ryu, Ernest K.},
	booktitle = {International Conference on Machine Learning},
	pages = {12098--12109},
	year = {2021},
}

@inproceedings{cai2022stochastic,
	title = {Stochastic {H}alpern iteration with variance reduction for
	         stochastic monotone inclusions},
	author = {Cai, Xufeng and Song, Chaobing and Guzm{\'a}n, Crist{\'o}bal
	          and Diakonikolas, Jelena},
	booktitle = {Advances in Neural Information Processing Systems},
	volume = {35},
	pages = {24766--24779},
	year = {2022},
}

@article{lieder2021convergence,
	title = {On the convergence rate of the {H}alpern-iteration},
	author = {Lieder, Felix},
	journal = {Optimization Letters},
	volume = {15},
	number = {2},
	pages = {405--418},
	year = {2021},
	publisher = {Springer},
}

@article{bravo2024stochastic,
	title = {Stochastic {H}alpern iteration in normed spaces and
	         applications to reinforcement learning},
	author = {Bravo, Mario and Contreras, Juan Pablo},
	journal = {arXiv preprint arXiv:2403.12338},
	year = {2024},
}

@inproceedings{lee2025near,
	title = {Near-optimal sample complexity for {MDP}s via anchoring},
	author = {Lee, Jongmin and Bravo, Mario and Cominetti, Roberto},
	booktitle = {International Conference on Machine Learning},
	pages = {32907--32929},
	year = {2025},
}

@article{alacaoglu2025towards,
	author = {Alacaoglu, Ahmet and Malitsky, Yura and Wright, Stephen J.},
	title = {Towards weaker variance assumptions for stochastic optimization
	         },
	journal = {arXiv preprint arXiv:2504.09951},
	year = {2025},
}

@article{Ljung1977,
	author = {Ljung, Lennart},
	title = {Analysis of recursive stochastic algorithms},
	journal = {IEEE Transactions on Automatic Control},
	year = {1977},
	volume = {22},
	number = {4},
	pages = {551--575},
	publisher = {IEEE},
}

@book{KushnerYin2003,
	author = {Kushner, Harold J. and Yin, G. George},
	title = {Stochastic approximation and recursive algorithms and
	         applications},
	volume = {35},
	year = {2003},
	publisher = {Springer},
}

@book{Borkar2008,
	author = {Borkar, Vivek S.},
	title = {Stochastic approximation: a dynamical systems viewpoint},
	year = {2008},
	publisher = {Cambridge University Press},
}

@article{Benaim1996,
	author = {Michel Bena{\"{\i}}m},
	title = {A dynamical system approach to stochastic approximations},
	journal = {SIAM Journal on Control and Optimization},
	year = {1996},
	volume = {34},
	number = {2},
	pages = {437--472},
	publisher = {SIAM},
}

@incollection{Benaim1999,
	author = {Michel Bena{\"{\i}}m},
	title = {Dynamics of stochastic approximation algorithms},
	booktitle = {S{\'e}minaire de Probabilit{\'e}s XXXIII},
	volume = {1709},
	pages = {1--68},
	year = {1999},
	publisher = {Springer},
}

@article{BorkarMeyn2000,
	author = {Borkar, Vivek S. and Meyn, Sean P.},
	title = {The {O}.{D}.{E}. method for convergence of stochastic
	         approximation and reinforcement learning},
	journal = {SIAM Journal on Control and Optimization},
	year = {2000},
	volume = {38},
	number = {2},
	pages = {447--469},
	publisher = {SIAM},
}

@article{AndrieuMoulinesPriouret2005,
	author = {Christophe Andrieu and {\'E}ric Moulines and Pierre Priouret},
	title = {Stability of stochastic approximation under verifiable
	         conditions},
	journal = {SIAM Journal on Control and Optimization},
	year = {2005},
	volume = {44},
	number = {1},
	pages = {283--312},
	publisher = {SIAM},
}

@article{PerkinsLeslie2012,
	author = {Steven Perkins and David S. Leslie},
	title = {Asynchronous stochastic approximation with differential
	         inclusions},
	journal = {Stochastic Systems},
	year = {2012},
	volume = {2},
	number = {2},
	pages = {409--446},
	publisher = {INFORMS},
}

@book{Khalil2002,
	author = {Khalil, Hassan K.},
	title = {Nonlinear systems},
	year = {2002},
	publisher = {Prentice Hall},
}

@inproceedings{SuttonEtAl2009ICML,
	author = {Sutton, Richard S. and Maei, Hamid R. and Precup, Doina and
	          Bhatnagar, Shalabh and Silver, David and Szepesv{\'a}ri, Csaba
	          and Wiewiora, Eric},
	title = {Fast gradient-descent methods for temporal-difference learning
	         with linear function approximation},
	booktitle = {International Conference on Machine Learning},
	year = {2009},
}

@inproceedings{MaeiEtAl2009NIPS,
	author = {Maei, Hamid R. and Szepesv{\'a}ri, Csaba and Bhatnagar,
	          Shalabh and Precup, Doina and Silver, David and Sutton, Richard
	          S.},
	title = {Convergent temporal-difference learning with arbitrary smooth
	         function approximation},
	booktitle = {Advances in Neural Information Processing Systems},
	year = {2009},
}

@inproceedings{MaeiEtAl2010ICML,
	author = {Maei, Hamid R. and Szepesv{\'a}ri, Csaba and Bhatnagar,
	          Shalabh and Sutton, Richard S.},
	title = {Toward off-policy learning control with function approximation},
	booktitle = {International Conference on Machine Learning},
	year = {2010},
}

@article{AbounadiBertsekasBorkar2001,
	author = {Abounadi, Jalal and Bertsekas, Dimitri P. and Borkar, Vivek S.
	          },
	title = {Learning algorithms for {M}arkov decision processes with
	         average cost},
	journal = {SIAM Journal on Control and Optimization},
	year = {2001},
	volume = {40},
	number = {3},
	pages = {681--698},
	publisher = {SIAM},
}

@inproceedings{DevrajMeyn2017NIPS,
	author = {Devraj, Adithya M. and Meyn, Sean P.},
	title = {Zap {Q}-learning},
	booktitle = {Advances in Neural Information Processing Systems},
	volume = {30},
	year = {2017},
}

@inproceedings{ChenEtAl2019arXiv,
	author = {Chen, Shuhang and Devraj, Adithya M. and Lu, Fan and Busic,
	          Ana and Meyn, Sean P.},
	booktitle = {Advances in Neural Information Processing Systems},
	pages = {16879--16890},
	title = {Zap {Q}-learning with nonlinear function approximation},
	volume = {33},
	year = {2020},
}

@article{BhatnagarEtAl2009NAC,
	author = {Bhatnagar, Shalabh and Sutton, Richard S. and Ghavamzadeh,
	          Mohammad and Lee, Mark},
	title = {Natural actor--critic algorithms},
	journal = {Automatica},
	year = {2009},
	volume = {45},
	number = {11},
	pages = {2471--2482},
	publisher = {Elsevier},
}

@book{Chen2002StochasticApproximation,
	author = {Han-Fu Chen},
	title = {Stochastic approximation and its applications},
	year = {2002},
	volume = {64},
	publisher = {Springer},
}

@inproceedings{AsadiLittman2017AltSoftmax,
	author = {Asadi, Kavosh and Littman, Michael L.},
	title = {An alternative softmax operator for reinforcement learning},
	booktitle = {International Conference on Machine Learning},
	pages = {243--252},
	year = {2017},
}

@book{KushnerClark1978,
	author = {Kushner, Harold J. and Clark, Dean S.},
	title = {Stochastic approximation methods for constrained and
	         unconstrained systems},
	year = {1978},
	volume = {26},
	publisher = {Springer-Verlag},
}

@inproceedings{an2018pid,
	title = {A {PID} controller approach for stochastic optimization of deep
	         networks},
	author = {An, Wangpeng and Wang, Haoqian and Sun, Qingyun and Xu, Jun
	          and Dai, Qionghai and Zhang, Lei},
	booktitle = {Proceedings of the IEEE conference on computer vision and
	             pattern recognition},
	pages = {8522--8531},
	year = {2018},
}

@article{chen2024accelerated,
	title = {Accelerated optimization in deep learning with a
	         proportional-integral-derivative controller},
	author = {Chen, Song and Liu, Jiaxu and Wang, Pengkai and Xu, Chao and
	          Cai, Shengze and Chu, Jian},
	journal = {Nature Communications},
	volume = {15},
	number = {1},
	pages = {10263},
	year = {2024},
	publisher = {Nature Publishing Group},
}

@inproceedings{wang2016dueling,
	title = {Dueling network architectures for deep reinforcement learning},
	author = {Wang, Ziyu and Schaul, Tom and Hessel, Matteo and Hasselt,
	          Hado and Lanctot, Marc and Freitas, Nando},
	booktitle = {International Conference on Machine Learning},
	pages = {1995--2003},
	year = {2016},
}

@article{hager1989updating,
	title = {Updating the inverse of a matrix},
	author = {Hager, William W},
	journal = {SIAM review},
	volume = {31},
	number = {2},
	pages = {221--239},
	year = {1989},
	publisher = {SIAM},
}

@book{liberzon2003switching,
	title = {Switching in systems and control},
	author = {Liberzon, Daniel},
	volume = {190},
	year = {2003},
	publisher = {Springer},
}

@article{lee2020unified,
	title = {A unified switching system perspective and convergence analysis
	         of Q-learning algorithms},
	author = {Lee, Donghwan and He, Niao},
	journal = {Advances in neural information processing systems},
	volume = {33},
	pages = {15556--15567},
	year = {2020},
}
